\documentclass[11pt,leqno]{article}
\usepackage{amsfonts}
\usepackage{amsmath}
\usepackage{amssymb}
\usepackage{amsthm}
\usepackage{latexsym}
\usepackage{mathrsfs}
\usepackage{amsthm}
\usepackage{mathtools}
\usepackage{bm}
\usepackage{verbatim}
\usepackage{color}

\newtheorem{theorem}{\bf Theorem}[section]
\newtheorem{lemma}[theorem]{\bf Lemma}
\newtheorem{corollary}[theorem]{\bf Corollary}
\newtheorem{proposition}[theorem]{\bf Proposition}

\newtheorem{remark}[theorem]{\bf Remark}

\numberwithin{equation}{section}

\newcommand{\ds}{\partial_s}

\makeatletter
\newcommand{\opnorm}{\@ifstar\@opnorms\@opnorm}
\newcommand{\@opnorms}[1]{%
  \left|\mkern-1.5mu\left|\mkern-1.5mu\left|
   #1
  \right|\mkern-1.5mu\right|\mkern-1.5mu\right|
}
\newcommand{\@opnorm}[2][]{%
  \mathopen{#1|\mkern-1.5mu#1|\mkern-1.5mu#1|}
  #2
  \mathclose{#1|\mkern-1.5mu#1|\mkern-1.5mu#1|}
}

\begin{document}
\vspace*{0ex}
\begin{center}
{\Large\bf
Non-stationary motion of an inextensible closed string
}
\end{center}

\begin{center}
Tatsuo Iguchi and Masahiro Takayama\footnote{Corresponding author}
\end{center}

\bigskip

\begin{abstract}
We consider a non-stationary motion of an inextensible closed string. 
The motion of the string is governed by nonlinear and nonlocal hyperbolic equations. 
The nonlocality of the equations is caused by the unknown tension of the string and reflects the inextensibility of the string. 
We show that the initial value problem to the equations of motion is well-posed locally in time. 
In the proof, it is crucial to show a strict positivity of the tension when the string is in motion. 
To this end, we derive a sharp lower bound of the Green's function for a periodic boundary value problem. 
We also derive a blow-up criterion of smooth solutions, which states that any singularities on the velocity or on the tension of the string 
do not take place before a curvature singularity appears. 
\end{abstract}

\section{Introduction}\label{sec:intro}
We are concerned with a motion of an inextensible closed string of length $2\pi$. 
The string is described as a curve $\bm{x}(s,t)=(x_1(s,t),x_2(s,t),x_3(s,t))$ 
at time $t$ with an arc length parameter $s\in\mathbb{T}=\mathbb{R}/2\pi\mathbb{Z}$. 
Let $\tau(s,t)$ be a scalar tension of the string at the point $\bm{x}(s,t)$ at time $t$. 
Then, the motion of the string without any external forces is described by 
\begin{align}
\label{Eq}
& \ddot{\bm{x}}-(\tau\bm{x}')' = \bm{0} \quad\mbox{in}\quad \mathbb{T}\times(0,T), \\
\label{Constraint}
& |\bm{x}'|=1 \quad\qquad\quad \mbox{in}\quad \mathbb{T}\times(0,T), 
\end{align}
where $\dot{\bm{x}}$ and $\bm{x}'$ denote the derivatives of $\bm{x}$ with respect to $t$ and $s$, respectively. 
The first equation is the equation of motion and the second one expresses the inextensibility of the string. 
For a derivation of these equations, we refer, for example, to Reeken \cite{Reeken1977} and Yong \cite{Yong2006}. 
In this paper, we are going to consider the initial value problem to these equations under the initial conditions 
\begin{equation}
\label{IC}
(\bm{x},\dot{\bm{x}})|_{t=0}
=(\bm{x}_0^\mathrm{in},\bm{x}_1^\mathrm{in}) \quad\mbox{in}\quad \mathbb{T}.
\end{equation}
Note that for the problem \eqref{Eq}--\eqref{IC} to have a regular solution, 
the initial data $(\bm{x}_0^\mathrm{in},\bm{x}_1^\mathrm{in})$ are forced to satisfy the conditions 
\begin{equation}
\label{Compatibility_Condition}
|\bm{x}_0^{\mathrm{in}\,\prime}|=1, \quad
\bm{x}_0^{\mathrm{in\,\prime}}\cdot\bm{x}_1^{\mathrm{in}\,\prime}=0 \quad\mbox{in}\quad \mathbb{T}. 
\end{equation}

The unknown quantities in this problem are not only the position vector $\bm{x}$ but also the tension $\tau$ of the string. 
Apparently, it is not quite clear how one can solve the tension $\tau$ from the equations \eqref{Eq} and \eqref{Constraint}. 
However, as in Preston \cite[Section 2.1]{Preston2011} and \c{S}eng\"{u}l and Vorotnikov \cite[Section 2.4]{SengulVorotnikov2017}, 
it is easily seen that for any regular solution $(\bm{x},\tau)$ of \eqref{Eq} and \eqref{Constraint} one has 
\begin{equation}\label{BVP}
-\tau''+|\bm{x}''|^2\tau = |\dot{\bm{x}}'|^2 \quad \mbox{in}\quad \mathbb{T}\times(0,T).
\end{equation}
Therefore, once the curvature vector $\bm{x}''(\cdot,t)$ and the velocity vector $\dot{\bm{x}}(\cdot,t)$ are given at time $t$, 
the tension $\tau(\cdot,t)$ at time $t$ is determined as a unique solution of the above ordinary differential equation. 
Plugging the solution $\tau$ into the equation of motion \eqref{Eq}, we obtain closed equations for the position vector $\bm{x}$, 
which are nonlinear and nonlocal hyperbolic equations.

As can be seen easily, for any regular solution $(\bm{x},\tau)$ of \eqref{Eq} the momentum is conserved, that is, 
$\int_{\mathbb{T}}\dot{\bm{x}}(s,t)\mathrm{d}s \equiv \int_{\mathbb{T}}\dot{\bm{x}}(s,0)\mathrm{d}s$. 
Moreover, for any regular solution $(\bm{x},\tau)$ of \eqref{Eq}--\eqref{Constraint} the kinetic energy is also conserved, that is, 
$\|\dot{\bm{x}}(t)\|_{L^2}\equiv\|\dot{\bm{x}}(0)\|_{L^2}$; see Lemma \ref{lem:CKE}. 
Now, let us assume that the initial velocity is constant, that is, $\bm{x}_1^\mathrm{in}(s)\equiv\bm{v}_0$ for some constant vector $\bm{v}_0$ 
and that $(\bm{x},\tau)$ is a regular solution of \eqref{Eq}--\eqref{IC}. 
Under the change of unknown variables $\bm{y}(s,t)=\bm{x}(s,t)-\bm{v}_0t$, the equations \eqref{Eq}--\eqref{Constraint} remain invariant, 
while the initial conditions \eqref{IC} are transformed into $(\bm{y},\dot{\bm{y}})|_{t=0}=(\bm{x}_0^\mathrm{in},\bm{0})$. 
Therefore, by the conservation of the kinetic energy we have $\bm{y}(s,t)\equiv \bm{x}_0^\mathrm{in}(s)$, 
so that $\bm{x}(s,t)=\bm{x}_0^\mathrm{in}(s)+\bm{v}_0t$. 
Plugging this into \eqref{BVP}, we obtain $\tau(s,t)\equiv0$. 
These arguments can be summarized as follows: If the initial velocity is a constant vector $\bm{v}_0$, 
then the initial value problem \eqref{Eq}--\eqref{IC} has a unique time global solution $(\bm{x},\tau)=(\bm{x}_0^{\mathrm{in}}(s)+\bm{v}_0t,0)$, 
which indicates that the string undergoes uniform linear motion. 
In view of this, it is sufficient to consider the case where the initial velocity $\bm{x}_1^\mathrm{in}(s)$ is not constant, 
that is, the case $\|\bm{x}_1^{\mathrm{in}\,\prime}\|_{L^2}>0$.

The following theorem is one of main results in this paper and ensures the well-posedness locally in time of the initial value problem \eqref{Eq}--\eqref{IC} 
in the $L^2$ Sobolev spaces $H^m$ of order $m$ in the case $m\geq3$.

\begin{theorem}\label{Th:LWP}
For any positive constant $M_*$ there exist a positive time $T$ and a positive constant $C$ such that if the initial data satisfy 
$(\bm{x}_0^{\mathrm{in}},\bm{x}_1^{\mathrm{in}})\in H^m\times H^{m-1}$ with an integer $m\geq3$, 
\eqref{Compatibility_Condition}, and 
\begin{equation}\label{Cond_ID} 
\|\bm{x}_1^{\mathrm{in}\,\prime}\|_{L^2}\geq M_*^{-1}, \quad 
\|\bm{x}_0^{\mathrm{in}}\|_{H^3}+\|\bm{x}_1^{\mathrm{in}}\|_{H^2}\leq M_*, 
\end{equation}
then the initial value problem \eqref{Eq}--\eqref{IC} has a unique solution $(\bm{x},\tau)$ on the time interval $[0,T]$ satisfying 
\begin{equation}\label{SolClass}
\bm{x}\in C^0([0,T];H^m)\cap C^1([0,T];H^{m-1}), \quad \tau \in C^0([0,T];H^m).
\end{equation}
Moreover, the tension satisfies $\tau(s,t)\geq C^{-1}$ for $(s,t)\in\mathbb{T}\times[0,T]$. 
\end{theorem}

\begin{remark}\label{Re:MainTh} \ 
{\rm (1)}
One may slightly generalize the equation of motion \eqref{Eq} by taking into account an effect of a constant gravity field $\bm{g}$ as 
\[
\ddot{\bm{x}}-(\tau\bm{x}')' = \bm{g} \quad\mbox{in}\quad \mathbb{T}\times(0,T).
\]
In fact, by a simple change of variables $\tilde{\bm{x}}(s,t)=\bm{x}(s,t)+\frac{1}{2}\bm{g}t^2$ from $\bm{x}$ to $\tilde{\bm{x}}$, 
one can reduce it to the case $\bm{g}=\bm{0}$ keeping the other equations unchanged. 

\smallskip
{\rm (2)}
Once a solution $(\bm{x},\tau)$ of \eqref{Eq}--\eqref{Constraint} satisfying the regularity \eqref{SolClass} is obtained, 
one can show further regularities of the solution in time by using the equations \eqref{Eq} and \eqref{BVP}, that is, 
$\bm{x}\in\bigcap_{j=0}^m C^j([0,T];H^{m-j})$ and $\tau\in\bigcap_{j=0}^{m-1} C^j([0,T];H^{m-j})$. 

\smallskip
{\rm (3)}
Our result does not include the critical case $m=2$ at the quasilinear regularity threshold. 
This is caused by the fact that one needs to control $\dot{\tau}$ the time derivative of the tension in the energy estimates. 
To evaluate it, we shall use the equation 
\[
-\dot{\tau}''+|\bm{x}''|^2\dot{\tau} = 2\dot{\bm{x}}'\cdot\ddot{\bm{x}}'-2\tau\bm{x}''\cdot\dot{\bm{x}}'' \quad \mbox{in}\quad \mathbb{T}\times(0,T),
\]
which is obtained by differentiating \eqref{BVP} with respect to $t$. 
In an evaluation of the right-hand side, we shall face to use of the norms $\|\ddot{\bm{x}}'(t)\|_{L^2}$ and $\|\dot{\bm{x}}''(t)\|_{L^2}$. 
This evaluation restricts us to the case $m\geq3$. 
\end{remark}

It is natural to expect that the solution of \eqref{Eq}--\eqref{IC} would generally have a singularity in finite time. 
It is interesting to clarify which kinds of singularities take place on the solution. 
The following theorem gives a blow-up criterion of smooth solutions and states that 
any singularities on the velocity or on the tension of the string do not take place before a curvature singularity appears. 
Note that $\bm{x}''$ is the curvature vector of the string.

\begin{theorem}\label{Th:BUC}
Let $(\bm{x},\tau)$ be a solution of \eqref{Eq}--\eqref{Constraint} in the class 
\[
\bm{x}\in C^0([0,T);H^m)\cap C^1([0,T);H^{m-1}), \quad \tau \in C^0([0,T);H^m)
\]
for some $T\in(0,\infty)$ with an integer $m\geq3$. 
If $T$ is maximal, that is, if $(\bm{x},\tau)$ cannot be continued to the solution in the above class beyond $T$, then we have 
\[
\int_0^T\|\bm{x}''(t)\|_{H^1}^2 \mathrm{d}t = \infty.
\]
\end{theorem}

In proving Theorems \ref{Th:LWP} and \ref{Th:BUC}, it is crucial to derive sharp lower and upper bounds of the tension $\tau$ of the string. 
Since the tension $\tau$ is determined as a unique solution of the ordinary differential equation \eqref{BVP}, we will consider the equation 
\begin{equation}\label{ODE}
-u''+au=f \quad\mbox{in}\quad \mathbb{T}
\end{equation}
for an unknown function $u$. 
For our applications, it is sufficient to consider the case where the coefficient $a\in L^1$ is non-negative and satisfies $\|a\|_{L^1}>0$. 
In such a case, the solution $u$ of \eqref{ODE} is unique so that the Green's function $G(s,\sigma)$ exists uniquely. 
Moreover, it is well-known that the Green's function is strictly positive. 
Although \eqref{ODE} is a simple and fundamental equation, to the best of the author's knowledge, 
there is no literature giving sharp quantitative estimates for lower and upper bounds of the Green's function in terms of $\|a\|_{L^1}$. 
The following theorem gives such estimates.

\begin{theorem}\label{Th:Green}
Suppose that $a\in L^1$ is non-negative and satisfies $\|a\|_{L^1}>0$. 
Let $G(s,\sigma)$ be the Green's function for the periodic boundary value problem \eqref{ODE}. 
Then, we have 
\begin{equation}\label{EstG}
\frac{1}{\|a\|_{L^1}\exp(\frac{\pi}{2}\|a\|_{L^1})} 
\leq G(s,\sigma) 
\leq  \frac{1+\frac{\pi}{2}\|a\|_{L^1}}{\|a\|_{L^1}}
\end{equation}
for any $(s,\sigma)\in\mathbb{T}\times\mathbb{T}$. 
\end{theorem}

\begin{remark}\label{Re:EstG} \ 
{\rm (1)} 
In the case $a(s)=\alpha^2$ with a positive constant $\alpha$, the Green's function can be written explicitly as 
$G(s,\sigma)=\frac{\cosh(\alpha(\pi+\sigma-s))}{2\alpha\sinh(\alpha\pi)}$ for $0\leq\sigma\leq s\leq2\pi$, so that we have 
$\min_{(s,t)\in\mathbb{T}\times\mathbb{T}} G(s,\sigma)=\frac{1}{2\alpha\sinh(\alpha\pi)}$. 
Therefore, the exponential decrease of the lower bound in \eqref{EstG} as $\|a\|_{L^1}\to\infty$ cannot be evitable in general. 

\smallskip
{\rm (2)} 
The upper bound in \eqref{EstG} is optimal. 
In fact, for any positive constant $\alpha$ there exists a sequence of non-negative functions $\{a_j\}$ satisfying $\|a_j\|_{L^1}=\alpha$ for any $j\in\mathbb{N}$ 
such that the corresponding sequence of Green's functions $\{G_j\}$ satisfies 
$\lim_{j\to\infty}\max_{(s,t)\in\mathbb{T}\times\mathbb{T}} G_j(s,\sigma)=\frac{1+\frac{\pi}{2}\alpha}{\alpha}$. 
For more details, we refer to Section \ref{sec:Gopt}. 
\end{remark}

Although \eqref{Eq}--\eqref{Constraint} are basic and fundamental equations for the motion of an inextensible string, 
there are few articles dealing with mathematical analysis of the equations. 
Particularly, there is no result in the setting of closed string. 
As for studies on the non-stationary motion in the setting of hanging string, 
where one end of the string is fixed and the other end is free, 
Reeken \cite{Reeken1979-1, Reeken1979-2} studied the motion of an inextensible string of infinite length 
in a gravity field under the setting that the string is hanging from a spatial infinity, that is, one end of the string is fixed at the point $(0,0,\infty)$.  
For technical reasons he assumed that the acceleration of a gravity vector $\bm{g}(s)\in C^\infty([0,\infty))$ is constant only for $s\in[0,l]$ and 
grows linearly beyond $s=l$ for some positive $l$, where the arc length parameter $s\in[0,\infty)$ is measured from the free end of the string. 
Under this non-physical condition, he proved the existence locally in time and uniqueness of the solution relying on the Nash--Moser implicit function theorem 
provided that the initial data are sufficiently close to a trivial stationary solution in some weighted Sobolev spaces 
with weights imposed at the free end of the string $s=0$. 
Later, Preston \cite{Preston2011} studied the motion of an inextensible string of finite length in a non-gravity field. 
Under this special situation, he proved the existence locally in time and uniqueness of the solution 
in some weighted Sobolev spaces without imposing any restrictions on the size of the initial data. 
However, he did not show the strong continuity of the solution in time. 
Recently, Iguchi and Takayama \cite{IguchiTakayama2024, IguchiTakayama2025, IguchiTakayama_arXiv:2502.16904} studied 
the motion of an inextensible string of finite length in a constant gravity field which is not necessarily zero. 
To address the well-posedness of the problem, they introduced a natural condition called {\it stability condition}. 
Under this stability condition on the initial data, they proved that the problem is well-posed locally in time in some weighted Sobolev spaces 
including the critical case at the quasilinear regularity threshold without imposing any restrictions on the size of the initial data, 
so that the strong continuity of the solution in time is also demonstrated. 
Note that the weights of the weighted Sobolev spaces used by Reeken \cite{Reeken1979-1, Reeken1979-2}, Preston \cite{Preston2011}, 
and Iguchi and Takayama \cite{IguchiTakayama2024, IguchiTakayama2025, IguchiTakayama_arXiv:2502.16904} may be apparently different, but are essentially the same.

We also mention a non-trivial stationary motion of an inextensible hanging string in a gravity field. 
In connection with the motion, there is a problem in determining the shape of the rotating string together with its angular velocity.  
A representative result on this problem was given by Kolodner \cite{Kolodner1955}, 
who proved that the corresponding nonlinear eigenvalue problem with a constant angular velocity $\omega$ has exactly $n$ non-trivial solutions 
if and only if $\omega$ satisfies $\omega_n<\omega\leq\omega_{n+1}$ with $\omega_n\equiv\sigma_n\sqrt{|\bm{g}|/4L}$, 
where $\sigma_n$ is the $n$-th zero of the Bessel function $J_0(z)$, $\bm{g}$ is the acceleration of gravity vector, and $L$ is the length of the string.

\medskip
The contents of this paper are as follows. 
In Section \ref{sec:Sturm-Liouville}, we derive an expression of the Green's function 
for the periodic boundary value problem \eqref{ODE} in terms of fundamental solutions. 
By choosing appropriately the fundamental solutions and evaluating them in detail, 
we develop sharp lower and upper bounds of the Green's function to prove Theorem \ref{Th:Green}. 
In Section \ref{sec:Estimates_P-BVP}, by applying results obtained in Section \ref{sec:Sturm-Liouville}, 
we derive estimates for the solution of \eqref{ODE}, which lead directly estimates for the tension $\tau$ including its sharp lower bound. 
In Section \ref{sec:Equivalence}, 
we show that under the necessary conditions \eqref{Compatibility_Condition} for the initial data, the solution $(\bm{x},\tau)$ of 
\eqref{Eq}, \eqref{BVP}, and \eqref{IC} satisfies necessarily the constraint \eqref{Constraint}, so that one can put \eqref{Constraint} 
behind in constructing the solution. 
In Section \ref{sec:WP}, we show that the transformed initial value problem \eqref{Eq}, \eqref{BVP}, and \eqref{IC} is well-posed locally in time 
in the Sobolev space and complete the proof of Theorem \ref{Th:LWP}. 
In Section \ref{sec:BUC}, we give a proof of Theorem \ref{Th:BUC}.

\medskip
\noindent
{\bf Notation}. \ 
For $1\leq p\leq\infty$, we denote by $L^p$ the Lebesgue space on $\mathbb{T}$. 
For an integer $m$ and $1\leq p\leq\infty$, we denote by $W^{m,p}$ the $L^p$ Sobolev space of order $m$ on $\mathbb{T}$ and put $H^m=W^{m,2}$. 
The norm of a Banach space $B$ is denoted by $\|\cdot\|_B$. 
The inner product in $L^2$ is denoted by $(\cdot,\cdot)_{L^2}$. 
$[P,Q]=PQ-QP$ denotes the commutator. 
We denote by $C(a_1, a_2, \ldots)$ a positive constant depending on $a_1, a_2, \ldots$. 
This constant $C(a_1, a_2, \ldots)$ may change from line to line. 
$f\lesssim g$ means that there exists a non-essential positive constant $C$ such that $f\leq Cg$ holds. 
$f\simeq g$ means that $f\lesssim g$ and $g\lesssim f$ hold. 
$a_1 \vee a_2 = \max\{a_1,a_2\}$.

\section{Green's function}\label{sec:Sturm-Liouville}
In this section, we consider the periodic boundary value problem \eqref{ODE}. 
As is well-known, the problem is equivalent to the equation 
\begin{equation}\label{SLP:Eq}
-u''+au=f \quad\mbox{in}\quad (0,2\pi)
\end{equation}
under the periodic boundary condition 
\begin{equation}\label{SLP:BC}
u(0)=u(2\pi), \quad u'(0)=u'(2\pi).
\end{equation}
Throughout this section, we assume that the coefficient $a\in L^1$ is non-negative and satisfies $\|a\|_{L^1}>0$. 
Then, the solution $u$ of \eqref{SLP:Eq}--\eqref{SLP:BC} is unique so that the corresponding Green's function exists uniquely, 
which will be denoted by $G(s,\sigma)$ in the following.

\subsection{Expression of the Green's function}\label{sec:GF}
We consider the corresponding homogeneous equation to \eqref{SLP:Eq}, that is, 
\begin{equation}\label{SLP:hEq}
-\varphi''+a\varphi=0 \quad\mbox{in}\quad (0,2\pi).
\end{equation}
Let $\bm{\varphi}=(\varphi_1,\varphi_2)^\mathrm{T}$ be fundamental solutions of this equation and denote their Wronskian by 
$W=\varphi_1\varphi_2'-\varphi_2\varphi_1'$, which is a non-zero constant. 
Then, any solution $u$ of \eqref{SLP:Eq} as well as its derivative $u'$ can be written as 
\begin{equation}\label{GSol}
\begin{pmatrix} u(s) \\ u'(s) \end{pmatrix}
= \Phi(s) \left\{ \begin{pmatrix} C_1 \\ C_2 \end{pmatrix} + A\int_0^s\bm{\varphi}(\sigma)f(\sigma)\mathrm{d}\sigma \right\}
\end{equation}
with constants $C_1$ and $C_2$, where 
\begin{equation}\label{PhiA}
\Phi(s)=
\begin{pmatrix}
 \varphi_1(s) & \varphi_2(s) \\
 \varphi_1'(s) & \varphi_2'(s)
\end{pmatrix}, \quad
A=\frac{1}{W}
\begin{pmatrix}
 0 & 1 \\
 -1 & 0
\end{pmatrix}.
\end{equation}
This solution satisfies the periodic boundary condition \eqref{SLP:BC} if and only if the constants $C_1$ and $C_2$ satisfy 
\[
(\Phi(0)-\Phi(2\pi))\begin{pmatrix} C_1 \\ C_2 \end{pmatrix}
= \Phi(2\pi)A\int_0^{2\pi}\bm{\varphi}(\sigma)f(\sigma)\mathrm{d}\sigma.
\]
Since $\varphi_1$ and $\varphi_2$ are linearly independent, 
the matrix $\Phi(0)-\Phi(2\pi)$ is invertible so that the constants $C_1$ and $C_2$ are uniquely determined. 
Plugging the constants into \eqref{GSol}, we obtain 
\begin{align*}
\begin{pmatrix} u(s) \\ u'(s) \end{pmatrix}
=
\Phi(s)(\Phi(0)-\Phi(2\pi))^{-1}\left\{ \Phi(0)A\int_0^s\bm{\varphi}(\sigma)f(\sigma)\mathrm{d}\sigma
 + \Phi(2\pi)A\int_s^{2\pi}\bm{\varphi}(\sigma)f(\sigma)\mathrm{d}\sigma \right\}.
\end{align*}
Here, we see that $(\Phi(0)-\Phi(2\pi))^{-1}\Phi(0)A = ((\Phi(0)-\Phi(2\pi))^{-1}\Phi(2\pi)A)^\mathrm{T}$. 
Therefore, we obtain 
\[
u(s) = \int_0^s\bm{\varphi}(s)\cdot B\bm{\varphi}(\sigma)f(\sigma)\mathrm{d}\sigma
 + \int_s^{2\pi}\bm{\varphi}(s)\cdot B^\mathrm{T}\bm{\varphi}(\sigma)f(\sigma)\mathrm{d}\sigma,
\]
where 
\begin{equation}\label{MatrixB}
B=(\Phi(0)-\Phi(2\pi))^{-1}\Phi(0)A.
\end{equation}
Consequently, we obtain an expression of the Green's function as 
\begin{equation}\label{ExpG}
G(s,\sigma) = 
\begin{cases}
 \bm{\varphi}(s)\cdot B\bm{\varphi}(\sigma) &\mbox{for}\quad 0\leq\sigma\leq s\leq2\pi, \\
 \bm{\varphi}(\sigma)\cdot B\bm{\varphi}(s) &\mbox{for}\quad 0\leq s\leq\sigma\leq2\pi.
\end{cases}
\end{equation}

Obviously, the Green's function is symmetric, that is, 
\begin{equation}
\label{symmetry_GF}
G(s,\sigma)=G(\sigma,s) \quad\mbox{for}\quad s,\sigma\in[0,2\pi].
\end{equation}
Since $\Phi(s)A\bm{\varphi}(s)=-\mathbf{e}_2$ and $\Phi(s)A\bm{\varphi}'(s)=\mathbf{e}_1$ hold with the standard basis $\{\mathbf{e}_1,\mathbf{e}_2\}$ 
of $\mathbb{R}^2$, we have $B\bm{\varphi}(0)=B^\mathrm{T}\bm{\varphi}(2\pi)$ and $B\bm{\varphi}'(0)=B^\mathrm{T}\bm{\varphi}'(2\pi)$. 
Therefore, the Green's function is periodic with respect to $s$ and $\sigma$, that is, 
\begin{equation}
\label{periodicity1_GF}
 \begin{cases}
  G(0,\sigma)=G(2\pi,\sigma) &\mbox{for}\quad \sigma\in[0,2\pi], \\
  G(s,0)=G(s,2\pi) &\mbox{for}\quad s\in[0,2\pi],
 \end{cases}
\end{equation}
and 
\begin{equation}
\label{periodicity2_GF}
 \begin{cases}
  (\partial_sG)(0,\sigma)=(\partial_sG)(2\pi,\sigma) &\mbox{for}\quad \sigma\in(0,2\pi), \\
  (\partial_\sigma G)(s,0)=(\partial_\sigma G)(s,2\pi) &\mbox{for}\quad s\in(0,2\pi).
 \end{cases}
\end{equation}
Moreover, we have the identity $B-B^\mathrm{T}=A$ so that $\bm{\varphi}'(s)\cdot(B-B^\mathrm{T})\bm{\varphi}(s)\equiv -1$, which implies 
the jump condition 
\begin{equation}
\label{jump_GF}
 \begin{cases}
  (\partial_sG)(\sigma+0,\sigma)-(\partial_sG)(\sigma-0,\sigma)=-1 &\mbox{for}\quad \sigma\in(0,2\pi), \\
  (\partial_\sigma G)(s,s-0)-(\partial_\sigma G)(s,s+0)=1 &\mbox{for}\quad s\in(0,2\pi).
 \end{cases}
\end{equation}

In evaluating the Green's function, we will use the expression \eqref{ExpG}. 
Although the Green's function is unique, there are many freedoms to choose the fundamental solutions $\bm{\varphi}$ and a better choice would lead 
sharper lower and upper bounds of the Green's function. 
In this paper, we choose them as solutions of the initial value problem to \eqref{SLP:hEq} under the initial conditions 
\begin{equation}\label{ICs}
\varphi_1(0)=0, \quad \varphi_1'(0)=1; \qquad \varphi_2(2\pi)=0, \quad \varphi_2'(2\pi)=-1.
\end{equation}
Here, we need to check that $\varphi_1$ and $\varphi_2$ are linearly independent. 
By the same argument as the proof of \cite[Lemma 3.1]{IguchiTakayama2024}, we have 
\begin{equation}\label{Est_phi}
\begin{cases}
 1\leq\varphi_1'(s) \leq \exp(\int_0^s \sigma a(\sigma)\mathrm{d}\sigma), \\
 s\leq\varphi_1(s) \leq s\exp(\int_0^s \sigma a(\sigma)\mathrm{d}\sigma)
\end{cases}
\end{equation}
for $s\in[0,2\pi]$. 
To evaluate $\varphi_2$, we put $\tilde{a}(s)=a(2\pi-s)$ and $\tilde{\varphi}_2(s)=\varphi_2(2\pi-s)$ for $s\in[0,2\pi]$.  
Then, we see that $\tilde{\varphi}_2$ is a unique solution of \eqref{SLP:hEq} with $a$ replaced by $\tilde{a}$ and satisfies the same initial conditions 
as those for $\varphi_1$. 
Moreover, we have $\int_0^{2\pi-s} \sigma \tilde{a}(\sigma)\mathrm{d}\sigma=\int_s^{2\pi} (2\pi-\sigma) a(\sigma)\mathrm{d}\sigma$. 
Therefore, we also have 
\begin{equation}\label{Est_psi}
\begin{cases}
 1\leq -\varphi_2'(s) \leq \exp(\int_s^{2\pi} (2\pi-\sigma) a(\sigma)\mathrm{d}\sigma), \\
 2\pi-s\leq\varphi_2(s) \leq (2\pi-s)\exp(\int_s^{2\pi} (2\pi-\sigma) a(\sigma)\mathrm{d}\sigma)
\end{cases}
\end{equation}
for $s\in[0,2\pi]$. 
Particularly, we see that the Wronskian $W$ satisfies 
\begin{equation}\label{Wronskian}
W = -\varphi_1(2\pi) = -\varphi_2(0) \leq -2\pi,
\end{equation}
so that they are in fact fundamental solutions of \eqref{SLP:hEq}. 
We put 
\begin{equation}\label{D}
D=\det(\Phi(0)-\Phi(2\pi))=W(2+\varphi_2'(0)-\varphi_1'(2\pi)),
\end{equation}
which is necessarily non-zero. 
Then, the matrix $B$ defined in \eqref{MatrixB} can be written as 
\[
B=\frac{1}{D}\begin{pmatrix} 1 & 1 \\ 1 & 1 \end{pmatrix} - \frac{1}{W}\begin{pmatrix} 0 & 0 \\ 1 & 0 \end{pmatrix}.
\]
Consequently, we get 
\begin{equation}\label{expression_GF}
G(s,\sigma) = \frac{(\varphi_1(s)+\varphi_2(s))(\varphi_1(\sigma)+\varphi_2(\sigma))}{W(2+\varphi_2'(0)-\varphi_1'(2\pi))}
 - 
\begin{cases}
 \dfrac{\varphi_2(s)\varphi_1(\sigma)}{W} &\mbox{for}\quad 0\leq\sigma\leq s\leq 2\pi, \\[2ex]
 \dfrac{\varphi_1(s)\varphi_2(\sigma)}{W} &\mbox{for}\quad 0\leq s\leq\sigma\leq 2\pi. 
\end{cases}
\end{equation}

\subsection{Quantitative estimates of the Green's function}
We proceed to prove Theorem \ref{Th:Green}, namely, we derive sharp lower and upper bounds of the Green's function relaying on its expression \eqref{expression_GF}. 
We first derive the upper bound.

\begin{lemma}
\label{lem:Est_phi+psi}
It holds that $\min_{s\in[0,2\pi]}(\varphi_1(s)+\varphi_2(s))\geq 2\pi$ and $\max_{s\in[0,2\pi]}(\varphi_1(s)+\varphi_2(s))=|W|$.  
\end{lemma}

\begin{proof}
Let us set $F_1(s)=\varphi_1(s)+\varphi_2(s)$ for $s\in[0,2\pi]$. 
The lower bound $F_1(s)\geq 2\pi$ follows easily from \eqref{Est_phi}--\eqref{Est_psi}. 
Since $F_1''(s)=a(s)F_1(s)\geq 0$ for $s\in(0,2\pi)$, we have $F_1(s)\leq\max\{F_1(0),F_1(2\pi)\}=|W|$ for $s\in[0,2\pi]$, 
where we used \eqref{ICs} and \eqref{Wronskian} in the last equality. 
\end{proof}

\begin{lemma}\label{lem:denominator_FS}
It holds that $2\pi\|a\|_{L^1} \leq \varphi_1'(2\pi)-\varphi_2'(0)-2 \leq |W|\|a\|_{L^1}$. 
\end{lemma}

\begin{proof}
Since $\varphi_1$ and $\varphi_2$ are solutions of \eqref{SLP:hEq} satisfying the initial conditions \eqref{ICs}, 
\begin{align*}
\varphi_1'(2\pi)-\varphi_2'(0)-2
&= \varphi_1'(2\pi)-\varphi_1'(0)+\varphi_2'(2\pi)-\varphi_2'(0) \\
&= \int_0^{2\pi}(\varphi_1''(\sigma)+\varphi_2''(\sigma))\mathrm{d}\sigma \\
&= \int_0^{2\pi}a(\sigma)(\varphi_1(\sigma)+\varphi_2(\sigma))\mathrm{d}\sigma. 
\end{align*}
Therefore, the desired estimate follows from Lemma \ref{lem:Est_phi+psi}. 
\end{proof}

\begin{lemma}\label{lem:key_inequality_FS}
It holds that $|W| \leq \frac{\pi}{2}(\varphi_1'(2\pi)-\varphi_2'(0)-2)+2\pi$. 
\end{lemma}

\begin{proof}
We set $\alpha_1=\varphi_1'(2\pi)-1$ and $\alpha_2=-\varphi_2'(0)-1$. 
As in the proof of Lemma \ref{lem:denominator_FS}, we have $\alpha_j=\int_0^{2\pi}a(s)\varphi_j(s)\mathrm{d}s$ for $j=1,2$, 
which together with \eqref{Est_phi}--\eqref{Est_psi} implies that $\alpha_1$ and $\alpha_2$ are both positive. 
Putting $\lambda_j=\frac{\alpha_j}{\alpha_1+\alpha_2}$ for $j=1,2$, we set $F_2(s)=\lambda_1\varphi_1(s)-\lambda_2\varphi_2(s)$ for $s\in[0,2\pi]$. 
Here, we note that $\lambda_1$ and $\lambda_2$ are chosen so that $\lambda_1+\lambda_2=1$ and $F_2'(0)=F_2'(2\pi)$ holds. 
By \eqref{Est_phi}--\eqref{Est_psi}, we have $F_2'(s)\geq1$ for $s\in[0,2\pi]$, 
$F_2(0)=-\lambda_2\varphi_2(0)<0$, and $F_2(2\pi)=\lambda_1\varphi_1(2\pi)>0$.   
Therefore, there exists an $s_*\in(0,2\pi)$ such that $F_2(s_*)=0$ 
and that $F_2''(s)=a(s)F_2(s)\leq 0$ for $0<s<s_*$ and $F_2''(s)=a(s)F_2(s)\geq 0$ for $s_*<s<2\pi$. 
Hence, we get $F_2'(s)\leq \max\{F_2'(0),F_2'(2\pi)\}$ for $s\in[0,2\pi]$. 
Now, noting that $|W|=\varphi_1(2\pi)=\varphi_2(0)$ we obtain 
\begin{align*}
|W|
&= F_2(2\pi)-F_2(0) \\
&= \int_0^{2\pi}F_2'(s)\mathrm{d}s \\
&\leq 2\pi( \lambda_1 F_2'(0)+\lambda_2 F_2'(2\pi) ) \\
&= 2\pi\lambda_1\lambda_2(\varphi_1'(2\pi)-\varphi_2'(0)-2)+2\pi \\
& \leq\frac{\pi}{2}(\varphi_1'(2\pi)-\varphi_2'(0)-2)+2\pi, 
\end{align*}
where we also used $\lambda_1\lambda_2\leq\frac{1}{4}$, which comes from $\lambda_1+\lambda_2=1$. 
This gives the desired estimate. 
\end{proof}

\begin{proposition}
\label{prop:GF_maximum}
It holds that $G(s,\sigma)  \leq \frac{1+\frac{\pi}{2}\|a\|_{L^1}}{\|a\|_{L^1}}$ for $(s,\sigma)\in[0,2\pi]\times[0,2\pi]$. 
\end{proposition}

\begin{proof}
Let $(s_0,\sigma_0)$ be a maximum point of the Green's function $G(s,\sigma)$. 
By the symmetry \eqref{symmetry_GF}, we may assume that $0\leq\sigma_0\leq s_0\leq 2\pi$.

We first consider the case $\sigma_0=0$. 
It follows from the expression \eqref{expression_GF} together with \eqref{ICs} and \eqref{Wronskian} that 
\begin{equation}\label{G(s,0)}
G(s,0) = \frac{\varphi_1(s)+\varphi_2(s)}{\varphi_1'(2\pi)-\varphi_2'(0)-2}. 
\end{equation}
Therefore, 
by Lemmas \ref{lem:Est_phi+psi}--\ref{lem:key_inequality_FS}, we see that  
\begin{align*}
G(s_0,0)
&= \frac{|W|}{\varphi_1'(2\pi)-\varphi_2'(0)-2} \\
&\leq \frac{\pi}{2}+\frac{2\pi}{\varphi_1'(2\pi)-\varphi_2'(0)-2} \\
&\leq \frac{\pi}{2}+\frac{1}{\|a\|_{L^1}}, 
\end{align*}
which gives the desired estimate.

We then consider the case $\sigma_0\in(0,2\pi]$. 
We extend the coefficient $a$ as a $2\pi$-periodic function, put $\tilde{a}(s)=a(s+\sigma_0)$, 
and denote the corresponding Green's function to $\tilde{a}$ by $\tilde{G}(s,\sigma)$. 
We also note that the Green's function $G(s,\sigma)$ can be extended naturally as a $2\pi$-periodic function with respect to $s$ and $\sigma$ 
thanks to the periodicity \eqref{periodicity1_GF} and \eqref{periodicity2_GF}. 
Then, by the uniqueness of the Green's function, we obtain $G(s,\sigma)=\tilde{G}(s-\sigma_0,\sigma-\sigma_0)$. 
Particularly, we get $G(s_0,\sigma_0)=\tilde{G}(s_0-\sigma_0,0)$. 
Since $\|\tilde{a}\|_{L^1}=\|a\|_{L^1}$, 
the desired estimate follows from the first case. 
\end{proof}

We proceed to derive the lower bound of the Green's function.

\begin{lemma}\label{lem:Min_phi+psi}
The function $\varphi_1(s)+\varphi_2(s)$ does not attain its minimum at $s=0$ nor $s=2\pi$. 
Moreover, it holds that $\min_{s\in[0,2\pi]}(\varphi_1(s)+\varphi_2(s)) \geq |W|\exp\left(-\frac{\pi}{2}\|a\|_{L^1}\right)$. 
\end{lemma}

\begin{proof}
Let us set $F_1(s)=\varphi_1(s)+\varphi_2(s)$ for $s\in[0,2\pi]$, and $s_0$ be a minimum point of $F_1(s)$. 
Suppose that $s_0=0$ or $s_0=2\pi$. 
In this case, we have $\min_{s\in[0,2\pi]}F_1(s)=|W|$, where we used \eqref{Wronskian}. 
Then, by Lemma \ref{lem:Est_phi+psi}, we get $F_1(s)\equiv |W|$, so that $\int_0^{2\pi}a(s)F_1(s)\mathrm{d}s=\int_0^{2\pi}F_1''(s)\mathrm{d}s=0$. 
On the other hand, Lemma \ref{lem:Est_phi+psi} implies that $\int_0^{2\pi}a(s)F_1(s)\mathrm{d}s \geq 2\pi\int_0^{2\pi}a(s)\mathrm{d}s=2\pi\|a\|_{L^1}$, 
which contradicts our assumption on the coefficient $a$. 
Therefore, we have $s_0\neq 0,2\pi$. 
This implies $F_1'(s_0)=0$, so that $\varphi_1'(s_0)=-\varphi_2'(s_0)$. 
Thus, evaluating the Wronskian at $s=s_0$ we get $|W|=F_1(s_0)\varphi_1'(s_0)$, which gives $F_1(s_0)=|W|/\varphi_1'(s_0)$. 
Putting $\lambda_1=\frac{2\pi-s_0}{2\pi}$ and $\lambda_2= \frac{s_0}{2\pi}$, we see by \eqref{Est_phi}--\eqref{Est_psi} that 
\begin{align*}
\varphi_1'(s_0)
&= \varphi_1'(s_0)^{\lambda_1}(-\varphi_2'(s_0))^{\lambda_2} \\
&\leq \exp\left(\lambda_1 s_0\int_0^{s_0}a(\sigma)\mathrm{d}\sigma\right)
 \exp\left(\lambda_2(2\pi-s_0)\int_{s_0}^{2\pi}a(\sigma)\mathrm{d}\sigma\right) \\
&= \exp\left(2\pi\lambda_1\lambda_2\|a\|_{L^1}\right) \\
& \leq\exp\left(\tfrac{\pi}{2}\|a\|_{L^1}\right), 
\end{align*}
where we also used $\lambda_1\lambda_2\leq\frac{1}{4}$, which comes from $\lambda_1+\lambda_2=1$. 
Therefore, we obtain the desired estimate. 
\end{proof}

\begin{proposition}\label{lem:GF_minimum}
It holds that $G(s,\sigma) \geq \frac{1}{\|a\|_{L^1}\exp(\frac{\pi}{2}\|a\|_{L^1})}$ 
for $(s,\sigma)\in[0,2\pi]\times[0,2\pi]$. 
\end{proposition}

\begin{proof}
Let $(s_0,\sigma_0)$ be a minimum point of the Green's function $G(s,\sigma)$. 
By the symmetry \eqref{symmetry_GF}, we may assume that $0\leq\sigma_0\leq s_0\leq 2\pi$. 
Moreover, by the same argument as in the proof of Proposition \ref{prop:GF_maximum} we may assume that $\sigma_0=0$. 
Then, in view of \eqref{G(s,0)} together with Lemmas \ref{lem:Min_phi+psi} and \ref{lem:denominator_FS}, we have 
\begin{align*}
G(s_0,0)
&= \frac{\varphi_1(s_0)+\varphi_2(s_0)}{\varphi_1'(2\pi)-\varphi_2'(0)-2} \\
&\geq \frac{1}{\|a\|_{L^1}\exp(\frac{\pi}{2}\|a\|_{L^1})}, 
\end{align*}
which gives the desired estimate. 
\end{proof}

We have now completed the proof of Theorem \ref{Th:Green}. 
Here, we also derive estimates for derivatives of the Green's function.

\begin{lemma}\label{lem:integral=1_GF}
It holds that $\int_0^{2\pi}G(s,\sigma)a(\sigma)\mathrm{d}\sigma=1$ for $s\in[0,2\pi]$. 
\end{lemma}

\begin{proof}
Note that $u(s)\equiv 1$ is a unique solution of the periodic boundary value problem \eqref{SLP:Eq}--\eqref{SLP:BC} with $f$ replaced by $a$. 
Therefore, we obtain the desired identity. 
\end{proof}

\begin{proposition}\label{Prop:EstDG}
It holds that $|\partial_sG(s,\sigma)|\leq1$ and $|\partial_\sigma G(s,\sigma)|\leq1$ for $(s,\sigma)\in[0,2\pi]\times[0,2\pi]$ 
satisfying $s\ne\sigma$ and $(s,\sigma)\ne(0,2\pi),(2\pi,0)$. 
\end{proposition}

\begin{proof}
By the symmetry \eqref{symmetry_GF}, it is sufficient to show that $|\partial_sG(s_0,\sigma_0)|\leq 1$ 
for any $(s_0,\sigma_0)\in[0,2\pi]\times[0,2\pi]$ satisfying $s_0\neq\sigma_0$ and $(s_0,\sigma_0)\neq(0,2\pi),(2\pi,0)$. 
Moreover, by the same argument as in the proof of Proposition \ref{prop:GF_maximum} we may assume that $\sigma_0=0$. 
Let $s_*$ be a minimum point of $G(s,0)$. 
In view of \eqref{G(s,0)} and Lemma \ref{lem:Min_phi+psi}, we have $s_*\neq 0,2\pi$, and hence $(\ds G)(s_*,0)=0$. 
Therefore, integrating $(\partial_s^2G)(s,0)=a(s)G(s,0)$, 
we have $(\partial_s G)(s_0,0)=\int_{s_*}^{s_0}a(s)G(s,0)\mathrm{d}s$, so that 
$|(\partial_s G)(s_0,0)|\leq\int_0^{2\pi}a(s)G(s,0)\mathrm{d}s=1$, 
where we used Lemma \ref{lem:integral=1_GF}. 
\end{proof}

\subsection{Optimality of the upper bound}\label{sec:Gopt}
We show that the upper bound of the Green's function given in Proposition \ref{prop:GF_maximum} is optimal. 
Let $\alpha$ be a fixed positive constant and $\delta$ a positive parameter satisfying $\delta^2<2\pi$, define a function $a(s)$ by 
\[
a(s)=
\begin{cases}
 \left(\frac{\alpha}{\delta}\right)^2 &\mbox{for}\quad 0<s<\delta^2, \\
 0 &\mbox{for}\quad \delta^2<s<2\pi,
\end{cases}
\]
and denote by $G(s,\sigma)=G(s,\sigma;\delta)$ the corresponding Green's function for the periodic boundary value problem \eqref{SLP:Eq}--\eqref{SLP:BC}. 
This function $a$ satisfies $\|a\|_{L^1}=\alpha^2$, which does not depend on $\delta$. 
The corresponding fundamental solutions $\varphi_1$ and $\varphi_2$ defined in Section \ref{sec:GF} are of the forms 
\begin{align*}
\varphi_1(s)
&= 
\begin{cases}
 \frac{\delta}{\alpha}\sinh\left(\frac{\alpha}{\delta} s\right) &\mbox{for}\quad 0\leq s\leq \delta^2, \\
 \frac{\delta}{\alpha}\sinh(\alpha\delta) + \cosh(\alpha\delta)(s-\delta^2) &\mbox{for}\quad \delta^2<s\leq 2\pi, 
\end{cases} \\
\varphi_2(s)
&= 
\begin{cases}
 (2\pi-\delta^2)\cosh\left(\frac{\alpha}{\delta}(\delta^2-s)\right)
 +\frac{\delta}{\alpha}\sinh\left(\frac{\alpha}{\delta} (\delta^2-s)\right) &\mbox{for}\quad 0\leq s\leq \delta^2, \\
 2\pi-s &\mbox{for}\quad \delta^2<s\leq 2\pi.
\end{cases}
\end{align*}
Therefore, in view of \eqref{expression_GF} we obtain 
\begin{align*}
\lim_{\delta\to+0}G(\pi,\pi;\delta)
&= \lim_{\delta\to+0} \left(\frac{(\varphi_1(\pi)+\varphi_2(\pi))^2}{W(2+\varphi_2'(0)-\varphi_1'(2\pi))}-\frac{\varphi_1(\pi)\varphi_2(\pi)}{W}\right) \\
&= \frac{1+\frac{\pi}{2}\alpha^2}{\alpha^2},
\end{align*}
so that 
\[
\lim_{\delta\to+0}\max_{(s,\sigma)\in[0,2\pi]\times[0,2\pi]}G(s,\sigma;\delta) 
=\frac{1+\frac{\pi}{2}\|a\|_{L^1}}{\|a\|_{L^1}}.
\]
This shows the criticality of the upper bound.

\section{Estimates for the tension $\tau$}\label{sec:Estimates_P-BVP}
We first consider the general periodic boundary value problem \eqref{ODE}, that is, 
\begin{equation}\label{ODE2}
-u''+au=f \quad\mbox{in}\quad \mathbb{T}.
\end{equation}
As in the previous section, we always assume that the coefficient $a\in L^1$ is non-negative and satisfies $\|a\|_{L^1}>0$.

\begin{lemma}\label{lem:BVP1}
For any $f\in L^1$, the periodic boundary value problem \eqref{ODE2} has a unique solution $u\in C^1(\mathbb{T})$, which satisfies 
$\|u\|_{L^\infty} \leq (\frac{\pi}{2}+1/\|a\|_{L^1})\|f\|_{L^1}$ and $\|u'\|_{L^\infty} \leq \|f\|_{L^1}$. 
If, in addition, $f$ is non-negative, then $u$ is also non-negative and satisfies 
$\min_{s\in\mathbb{T}}u(s) \geq \frac{1}{\|a\|_{L^1}\exp(\frac{\pi}{2}\|a\|_{L^1})}\|f\|_{L^1}$. 
\end{lemma}

\begin{proof}
The solution is expressed as $u(s)=\int_{\mathbb{T}}G(s,\sigma)f(\sigma)\mathrm{d}\sigma$, 
so that the first and the last estimates follow directly from Theorem \ref{Th:Green}. 
Since $u'(s)=\int_{\mathbb{T}}\ds G(s,\sigma)f(\sigma)\mathrm{d}\sigma$, the second one follows from Proposition \ref{Prop:EstDG}. 
\end{proof}

\begin{lemma}\label{lem:BVP2}
For any $f=F'$ with $F\in L^1$, the periodic boundary value problem \eqref{ODE2} has a unique solution $u\in W^{1,1}$, 
which satisfies $\|u\|_{L^\infty} \leq \|F\|_{L^1}$. 
If, in addition, $F\in L^2$, then $u\in H^1$ and satisfies $\|u'\|_{L^2} \leq \|F\|_{L^2}$. 
\end{lemma}

\begin{proof}
Since the solution is expressed as $u(s)=-\int_{\mathbb{T}}(\partial_\sigma G(s,\sigma))F(\sigma)\mathrm{d}\sigma$, 
the first estimate follows from Proposition \ref{Prop:EstDG}. 
By taking the $L^2$ inner product of the equation with $u$ and integration by parts, we see that 
$\|u'\|_{L^2}^2 + (au,u)_{L^2}=-(F,u')_{L^2} \leq \|F\|_{L^2}\|u'\|_{L^2}$, 
which implies the second one. 
\end{proof}

\begin{corollary}\label{Cor:BVP}
There exists an absolute positive constant $C$ such that for any $f\in H^{-1}$, 
the periodic boundary value problem \eqref{ODE2} has a unique solution $u\in H^1$, 
which satisfies $\|u\|_{L^\infty} \leq C(1+1/\|a\|_{L^1})\|f\|_{H^{-1}}$ and $\|u'\|_{L^2} \leq C\|f\|_{H^{-1}}$. 
\end{corollary}

\begin{proof}
It follows from Lemmas \ref{lem:BVP1} and \ref{lem:BVP2} that for any $f=f_1+f_2'$ with $f_1\in L^1$ and $f_2\in L^2$, 
\eqref{ODE2} has a unique solution $u\in H^1$, which satisfies 
$\|u\|_{L^\infty} \lesssim (1+1/\|a\|_{L^1})\|f_1\|_{L^1}+\|f_2\|_{L^2}$ and $\|u'\|_{L^2} \leq \|f_1\|_{L^1}+\|f_2\|_{L^2}$. 
These estimates together with the equivalence 
\begin{equation}\label{equi_norm}
\|f\|_{H^{-1}} \simeq \inf\{ \|f_1\|_{L^1}+\|f_2\|_{L^2} \,|\, f=f_1+f_2'\}
\end{equation}
imply the desired estimates. 
\end{proof}

We proceed to consider the equation \eqref{BVP} of the tension $\tau$, in which the time $t$ appears as a parameter. 
Therefore, in the following we fix the time $t$ arbitrarily, drop it in the notations, and consider the equation 
\begin{equation}\label{BVP2}
-\tau''+|\bm{x}''|^2\tau = |\dot{\bm{x}}'|^2 \quad\mbox{in}\quad \mathbb{T}.
\end{equation}
The following lemma follows directly from Lemmas \ref{lem:BVP1} and \ref{BVP2}.

\begin{lemma}\label{lem:tau1}
Suppose that $\bm{x}$ is given such that $\bm{x}'',\dot{\bm{x}}' \in L^2$ and $\|\bm{x}''\|_{L^2}>0$. 
Then, the tension $\tau$ is uniquely determined as a solution of \eqref{BVP2} and satisfies 
\[
\begin{cases}
 \displaystyle
 \|\tau\|_{L^\infty} \leq \left(\frac{\pi}{2}+\frac{1}{\|\bm{x}''\|_{L^2}^2}\right)\|\dot{\bm{x}}'\|_{L^2}^2, \quad
  \|\tau'\|_{L^\infty} \leq \|\dot{\bm{x}}'\|_{L^2}^2, \\[2.5ex]
 \displaystyle
 \min_{s\in\mathbb{T}}\tau(s) \geq \frac{\|\dot{\bm{x}}'\|_{L^2}^2}{\|\bm{x}''\|_{L^2}^2 \exp(\frac{\pi}{2}\|\bm{x}''\|_{L^2}^2)}.
\end{cases}
\]
\end{lemma}

\begin{remark}\label{Re:tau}
As we will see in Lemma \ref{lem:Below}, for any regular solution $(\bm{x},\tau)$ of \eqref{Eq}--\eqref{IC} one has 
$\|\dot{\bm{x}}'(t)\|_{L^2}\geq\|P_1\bm{x}_1^{\mathrm{in}}\|_{L^2}$, where $P_1$ is an orthogonal projection onto the non-constant component in $L^2$. 
Since we have assumed $\|\bm{x}_1^{\mathrm{in}\,\prime}\|_{L^2}>0$, one has $\|P_1\bm{x}_1^{\mathrm{in}}\|_{L^2}>0$. 
Hence, the lower bound in the above lemma ensures that the tension $\tau$ is strictly positive. 
\end{remark}

\begin{lemma}\label{lem:tau2}
For any integer $m\geq2$ and a positive constant $M$ there exists a positive constant $C(m,M)$ such that if $\bm{x}'',\dot{\bm{x}}' \in H^{(m-2)\vee1}$ 
satisfy $\|\bm{x}''\|_{L^2}\geq M^{-1}$ and $\|(\bm{x}'',\dot{\bm{x}}')\|_{H^1}\leq M$, then we have 
$\|\tau\|_{H^m} \leq C(m,M)\|(\bm{x}'',\dot{\bm{x}}')\|_{H^{m-2}}$. 
\end{lemma}

\begin{proof}
It follows from Lemma \ref{lem:tau1} that $\|\tau\|_{H^1} \lesssim \|\dot{\bm{x}}'\|_{L^2}$. 
For $k=2,3,\ldots,m$, by using the equation \eqref{BVP2} and standard calculus inequalities we see that 
\begin{align*}
\|\tau''\|_{H^{k-2}}
&\lesssim \|\tau\|_{H^{k-2}}\|\bm{x}''\|_{L^\infty}^2 + \|\tau\|_{L^\infty}\|\bm{x}''\|_{L^\infty}\|\bm{x}''\|_{H^{k-2}}
 + \|\dot{\bm{x}}'\|_{L^\infty}\|\dot{\bm{x}}'\|_{H^{k-2}} \\
&\lesssim \|\tau\|_{H^{k-2}} + \|(\bm{x}'',\dot{\bm{x}}')\|_{H^{k-2}},
\end{align*}
so that $\|\tau\|_{H^k} \lesssim \|\tau\|_{H^{k-2}} + \|(\bm{x}'',\dot{\bm{x}}')\|_{H^{k-2}}$. 
Using this inductively, we obtain the desired estimate. 
\end{proof}

In order to evaluate the time derivative $\dot{\tau}$ of the tension, we differentiate \eqref{BVP2} with respect to $t$ to obtain 
\begin{equation}\label{DtBVP}
-\dot{\tau}''+|\bm{x}''|^2\dot{\tau} = 2\dot{\bm{x}}'\cdot\ddot{\bm{x}}'-2\tau\bm{x}''\cdot\dot{\bm{x}}'' \quad \mbox{in}\quad \mathbb{T}.
\end{equation}
The following lemma follows from directly Lemma \ref{lem:BVP1}.

\begin{lemma}\label{lem:Dt_tau}
For any positive constant $M$ there exists a positive constant $C(M)$ such that if $\|\bm{x}''\|_{L^2}\geq M^{-1}$, then we have 
$\|\dot{\tau}\|_{W^{1,\infty}} 
\leq C(M)(\|\dot{\bm{x}}'\|_{L^2}\|\ddot{\bm{x}}'\|_{L^2}+\|\dot{\bm{x}}'\|_{L^2}^2\|\bm{x}''\|_{L^2}\|\dot{\bm{x}}''\|_{L^2})$. 
\end{lemma}

\section{Equivalence of the systems}\label{sec:Equivalence}
In constructing the solution $(\bm{x},\tau)$ of the problem \eqref{Eq}--\eqref{IC}, 
we will use \eqref{BVP} as the equation for $\tau$ and put the constraint \eqref{Constraint} behind. 
However, its solution does not necessarily satisfy \eqref{Constraint}. 
In this section, we will show that if the initial data satisfy the necessary conditions \eqref{Compatibility_Condition}, 
then the solution of \eqref{Eq}, \eqref{BVP}, and \eqref{IC} satisfies necessarily \eqref{Constraint}. 
The following lemma is a corollary of Fenchel's theorem.

\begin{lemma}\label{lem:Fenchel}
For any $\bm{x}\in H^2$ satisfying $|\bm{x}'(s)|\equiv 1$, 
we have $\|\bm{x}''\|_{L^2} \geq \sqrt{2\pi}$. 
\end{lemma}

\begin{proof}
Since $\bm{x}:\mathbb{T}\to\mathbb{R}^3$ is a closed curve parameterized by its arc length, the curvature $\kappa(s)$ is given by $\kappa(s)=|\bm{x}''(s)|$. 
Fenchel's theorem implies that $\int_{\mathbb{T}}\kappa(s)\mathrm{d}s\geq 2\pi$; see Fenchel \cite{Fenchel1929}. 
Therefore, we have $2\pi\leq\int_{\mathbb{T}}\kappa(s)\mathrm{d}s=\|\bm{x}''\|_{L^1}\leq\sqrt{2\pi}\|\bm{x}''\|_{L^2}$, which gives the desired estimate. 
\end{proof}

\begin{lemma}\label{lem:CKE}
For any regular solution $(\bm{x},\tau)$ of \eqref{Eq}--\eqref{Constraint}, we have $\|\dot{\bm{x}}(t)\|_{L^2}\equiv\|\dot{\bm{x}}(0)\|_{L^2}$. 
\end{lemma}

\begin{proof}
Differentiating \eqref{Constraint} with respect to $t$, we have $\dot{\bm{x}}'\cdot\bm{x}'=0$. 
By \eqref{Eq} and integration by parts, we see that 
$\frac12\frac{\mathrm{d}}{\mathrm{d}t}\|\dot{\bm{x}}\|_{L^2}^2=(\dot{\bm{x}},\ddot{\bm{x}})_{L^2}=(\dot{\bm{x}},(\tau\bm{x}')')_{L^2}
=-(\dot{\bm{x}}',\tau\bm{x}')_{L^2}=0$, which implies the desired conservation of the kinetic energy. 
\end{proof}

\begin{lemma}\label{lem:Below}
Let $P_0$ be an orthogonal projection in $L^2$ defined by $P_0u=\frac{1}{2\pi}\int_{\mathbb{T}}u(s)\mathrm{d}s$ and put $P_1=\mathrm{Id}-P_0$. 
For any regular solution $(\bm{x},\tau)$ of \eqref{Eq}--\eqref{IC}, we have $\|\dot{\bm{x}}'(t)\|_{L^2}\geq\|P_1\bm{x}_1^{\mathrm{in}}\|_{L^2}$. 
\end{lemma}

\begin{proof}
Putting $\bm{y}(s,t)=\bm{x}(s,t)-t(P_0\bm{x}_1^{\mathrm{in}})$, 
we see that $(\bm{y},\tau)$ is a solution of \eqref{Eq}--\eqref{Constraint} and that $\dot{\bm{y}}(\cdot,0)=P_1\bm{x}_1^{\mathrm{in}}$. 
Particularly, we have $\int_\mathbb{T}\dot{\bm{y}}(s,0)\mathrm{d}s=0$, which together with the conservation of the momentum implies 
$\int_\mathbb{T}\dot{\bm{y}}(s,t)\mathrm{d}s\equiv0$. 
Therefore, we can apply Poincar\'{e}'s inequality and the conservation of the kinetic energy ensured by Lemma \ref{lem:CKE} to obtain 
$\|\dot{\bm{x}}'(t)\|_{L^2}=\|\dot{\bm{y}}'(t)\|_{L^2}\geq\|\dot{\bm{y}}(t)\|_{L^2}=\|\dot{\bm{y}}(0)\|_{L^2}=\|P_1\bm{x}_1^{\mathrm{in}}\|_{L^2}$. 
\end{proof}

\begin{proposition}\label{Prop:Equivalence}
Let $(\bm{x},\tau)$ be a solution of the transformed problem \eqref{Eq}, \eqref{BVP}, and \eqref{IC} satisfying 
$\bm{x}\in L^\infty(0,T;H^3)$ and $\dot{\bm{x}}\in L^\infty(0,T;H^2)$. 
If the initial data $(\bm{x}_0^{\mathrm{in}},\bm{x}_1^{\mathrm{in}})$ satisfy \eqref{Compatibility_Condition} and $\|\bm{x}_1^{\mathrm{in}\,\prime}\|_{L^2}>0$, 
then we have $|\bm{x}'(s,t)|\equiv 1$ for $(s,t)\in\mathbb{T}\times[0,T]$. 
\end{proposition}

\begin{proof}
Let us define $T_*$ by 
\[
T_*=\sup\left\{T'\in[0,T] \;\Bigm|\; \min_{t\in[0,T']}\|\bm{x}''(t)\|_{L^2}>0 \ \mbox{ and } \min_{(s,t)\in\mathbb{T}\times[0,T']}\tau(s,t)>0\right\}. 
\]
Since $|\bm{x}'(s,0)|\equiv 1$, Lemma \ref{lem:Fenchel} implies $\|\bm{x}''(0)\|_{L^2}\geq\sqrt{2\pi}$. 
Moreover, by Lemma \ref{lem:tau1} we have $\tau(s,0)\gtrsim\|\dot{\bm{x}}'(0)\|_{L^2}^2 = \|\bm{x}_1^{\mathrm{in}\,\prime}\|_{L^2}^2>0$.  
Therefore, by the continuity of $\|\bm{x}''(t)\|_{L^2}$ and $\min_{s\in\mathbb{T}}\tau(s,t)$, we see that $0<T_*\leq T$.

We first show that 
\begin{equation}\label{Constraint_T*}
|\bm{x}'(s,t)|=1 \quad\mbox{for any} \quad (s,t)\in\mathbb{T}\times[0,T_*]. 
\end{equation}
Put $h(s,t)=|\bm{x}'(s,t)|^2-1$. 
By a straightforward calculation, we have
\[
\begin{cases}
 \ddot{h} = \tau h''+2\tau'h'+2\tau''h &\mbox{in} \quad \mathbb{T}\times(0,T'), \\
 (h,\dot{h})|_{t=0} = (0,0) &\mbox{in} \quad \mathbb{T},
\end{cases}
\]
where we used \eqref{Compatibility_Condition} to derive the homogeneous initial conditions. 
Therefore, putting $E(t)=\int_{\mathbb{T}}(|\dot{h}|^2+\tau|h'|^2+|h|^2)\mathrm{d}s$, we see that 
\[
\frac{\mathrm{d}}{\mathrm{d}t}E(t) = \int_{\mathbb{T}}(2\tau'\dot{h}h'+(4\tau''+2)\dot{h}h+\dot{\tau}|h'|^2)\mathrm{d}s. 
\]
Take $T'\in(0,T_*)$ arbitrarily. 
Then, we have $\|\bm{x}''(t)\|_{L^2}\gtrsim1$ and $\tau(s,t)\gtrsim1$ for $(s,t)\in\mathbb{T}\times[0,T']$, 
so that Lemmas \ref{lem:tau1}, \ref{lem:tau2}, and \ref{lem:Dt_tau} yield 
$\|\tau(t)\|_{W^{2,\infty}} \lesssim 1$ and $\|\dot{\tau}(t)\|_{L^\infty} \lesssim 1$ for $t\in[0,T']$. 
Therefore, we get $E(t)\simeq \|(h,\dot{h},h')(t)\|_{L^2}^2$ and that $\frac{\mathrm{d}}{\mathrm{d}t}E(t)\lesssim E(t)$ for $t\in[0,T']$, 
which together with the initial condition $E(0)=0$ and Gronwall's inequality yield $\|h(t)\|_{L^2}=0$ for $t\in[0,T']$. 
Since this equality holds for any $T'\in(0,T_*)$, by the continuity of $\|h(t)\|_{L^2}$ we get $\|h(t)\|_{L^2}=0$ for $t\in[0,T_*]$, 
which gives \eqref{Constraint_T*}.

We next show that $T_*=T$. 
Suppose that $T_*<T$. 
By \eqref{Constraint_T*}, we see that $(\bm{x},\tau)$ is a solution of \eqref{Eq}--\eqref{IC} on the time interval $[0,T_*]$, 
so that by Lemma \ref{lem:Fenchel}, we have $\|\bm{x}''(t)\|_{L^2}\geq \sqrt{2\pi}$ for $t\in[0,T_*]$. 
Moreover, by Lemmas \ref{lem:tau1} and \ref{lem:Below} we get 
$\min_{s\in\mathbb{T}}\tau(s,t)\gtrsim\|\dot{\bm{x}}'(t)\|_{L^2}^2 \geq \|P_1\bm{x}_1^{\mathrm{in}}\|_{L^2}^2>0$ for $t\in[0,T_*]$. 
Thus, by the continuity we have $\|\bm{x}''(t)\|_{L^2}\geq c$ and $\min_{s\in\mathbb{T}}\tau(s,t)\geq c$ for $t\in[0,T_*+\delta]$ 
with some small positive constants $\delta$ and $c$, which contradicts with the definition of $T_*$. 
Therefore, we get $T_*=T$. 
Combining this with \eqref{Constraint_T*} gives the desired result. 
\end{proof}

\section{Well-posedness of the problem}\label{sec:WP}
By taking Proposition \ref{Prop:Equivalence} into account, 
in order to prove Theorem \ref{Th:LWP} it is sufficient to show the well-posedness of the initial value problem 
\begin{align}\label{TFP}
 &\ddot{\bm{x}}-(\tau\bm{x}')' = \bm{0} \makebox[4.06em]{}\mbox{in}\quad \mathbb{T}\times(0,T), \\
  \label{TFBVP}
 &-\tau''+|\bm{x}''|^2\tau = |\dot{\bm{x}}'|^2    \quad\mbox{in}\quad \mathbb{T}\times(0,T)
\end{align}
under the initial conditions 
\begin{equation}\label{IC_bis}
(\bm{x},\dot{\bm{x}})|_{t=0}
=(\bm{x}_0^\mathrm{in},\bm{x}_1^\mathrm{in}) \quad\mbox{in}\quad \mathbb{T}.
\end{equation}

\begin{proposition}\label{Prop:WP}
For any constant $M_*\geq1$ there exist a positive time $T$ and a positive constant $C$ such that if the initial data satisfy 
$(\bm{x}_0^{\mathrm{in}},\bm{x}_1^{\mathrm{in}})\in H^m\times H^{m-1}$ with an integer $m\geq3$, 
$\|\bm{x}_0^{\mathrm{in}\,\prime\prime}\|_{L^2}\geq M_{*}^{-1}$, $\|\bm{x}_1^{\mathrm{in}\,\prime}\|_{L^2}>M_*^{-1}$, and 
$\|\bm{x}_0^{\mathrm{in}}\|_{H^3}+\|\bm{x}_1^{\mathrm{in}}\|_{H^2}\leq M_*$, 
then the initial value problem \eqref{TFP}--\eqref{IC_bis} has a unique solution $(\bm{x},\tau)$ on the time interval $[0,T]$ satisfying 
\begin{equation}\label{SolClass2}
\bm{x}\in C^0([0,T];H^m)\cap C^1([0,T];H^{m-1}), \quad \tau \in C^0([0,T];H^m).
\end{equation}
Moreover, the tension satisfies $\tau(s,t)\geq C^{-1}$ for $(s,t)\in\mathbb{T}\times[0,T]$. 
\end{proposition}

\begin{proof}
{\bf Construction of approximate solutions}. 
We are going to construct a sequence of approximate solutions $\{(\bm{x}_n,\tau_n)\}$. 
Put $\bm{x}_0(\cdot,t)=\bm{x}_0^\mathrm{in}+t\mathrm{e}^{-(1+|D|)t}\bm{x}_1^\mathrm{in}$, where $D=\frac{1}{\mathrm{i}}\frac{\partial}{\partial s}$ and 
we used the standard notation of Fourier multipliers. 
Then, $\bm{x}_0 \in \bigcap_{j=0}^m C^j([0,\infty);H^{m-j})$ and it satisfies the initial conditions \eqref{IC_bis} and 
\[
\begin{cases}
 \|\bm{x}_0(t)\|_{H^m}+\|\dot{\bm{x}}_0(t)\|_{H^{m-1}} \leq C(\|\bm{x}_0^\mathrm{in}\|_{H^m}+\|\bm{x}_1^\mathrm{in}\|_{H^{m-1}}), \\
 \|\bm{x}_0''(t)-\bm{x}_0^{\mathrm{in}\,\prime\prime}\|_{L^2} \leq CM_*t, \quad \|\dot{\bm{x}}_0'(t)-\bm{x}_1^{\mathrm{in}\,\prime}\|_{L^2} \leq CM_*t
\end{cases}
\]
for $t\geq0$ with an absolute positive constant $C$. 
Therefore, putting $M_0=2M_*$ and choosing $T_0>0$ so small that $2CM_*^2T_0\leq1$, 
we obtain $M_0^{-1} \leq \|\bm{x}_0''(t)\|_{L^2} \leq M_0$ and $M_0^{-1} \leq \|\dot{\bm{x}}_0'(t)\|_{L^2} \leq M_0$ for $t\in[0,T_0]$. 
Let $T\in(0,T_0]$ and suppose that $\bm{x}_n\in\bigcap_{j=0}^2C^j([0,T];H^{m-j})$ is given such that it satisfies the initial conditions \eqref{IC_bis} and 
\begin{equation}\label{SolEst}
 \begin{cases}
  M_0^{-1} \leq \|\bm{x}_n''(t)\|_{L^2} \leq M_0, \quad M_0^{-1} \leq \|\dot{\bm{x}}_n'(t)\|_{L^2} \leq M_0, \\
  \|(\dot{\bm{x}}_n,\bm{x}_n')(t)\|_{H^2} \leq M_1, \quad \|\ddot{\bm{x}}_n(t)\|_{H^1} \leq M_2 \quad\mbox{for}\quad t\in[0,T]
 \end{cases}
\end{equation}
with positive constants $M_1$ and $M_2$, which together with $T$ will be determined later. 
In view of Lemmas \ref{lem:tau1}, \ref{lem:tau2}, and \ref{lem:Dt_tau}, we define $\tau_n$ as a unique solution of the periodic boundary value problem 
\begin{equation}\label{BVPn}
-\tau_n''+|\bm{x}_n''|^2\tau_n=|\dot{\bm{x}}_n'|^2 \quad\mbox{in}\quad \mathbb{T}\times(0,T).
\end{equation}
Then, $\tau_n\in C^0([0,T];H^m)\cap C^1(\mathbb{T}\times[0,T])$ and it is strictly positive. 
Therefore, by a classical theory on the initial value problem to hyperbolic equations we can define $\bm{x}_{n+1}$ as a unique solution of the initial value problem 
\begin{equation}\label{IVPn}
 \begin{cases}
  \ddot{\bm{x}}_{n+1}=(\tau_n\bm{x}_{n+1}')' &\mbox{in}\quad \mathbb{T}\times(0,T), \\
  (\bm{x}_{n+1},\dot{\bm{x}}_{n+1})|_{t=0}=(\bm{x}_0^\mathrm{in},\bm{x}_1^\mathrm{in}) &\mbox{in}\quad \mathbb{T}
 \end{cases}
\end{equation}
in the class $\bm{x}_{n+1}\in\bigcap_{j=0}^2C^j([0,T];H^{m-j})$.

\medskip
{\bf Uniform bounds of the approximate solutions}. 
We proceed to show that $\bm{x}_{n+1}$ also satisfies \eqref{SolEst} with some constants $M_1, M_2$, and $T$ independent of $n$. 
In the following, we simply denote the constants $C_0=C(M_0)$, $C_1=C(M_0,M_1)$, $C_2=C(M_0,M_1,M_2)$, and $C_3=C(m,M_0,M_1,M_2)$. 
These constants may change from line to line. 
By Lemmas \ref{lem:tau1}, \ref{lem:tau2}, and \ref{lem:Dt_tau}, we have $\min_{s\in\mathbb{T}}\tau_n(s,t)\geq C_0^{-1}$, 
$\|\tau_n(t)\|_{W^{1,\infty}} \leq C_0$, $\|\tau_n(t)\|_{H^3} \leq C_1$, and $\|\dot{\tau}_n(t)\|_{W^{1,\infty}}\leq C_2$ for $t\in[0,T]$.

Putting $\bm{y}_{n+1}=\partial_s^j\bm{x}_{n+1}$ and differentiating \eqref{IVPn} $j$-times with respect to $s$, we obtain 
\begin{equation}\label{Eqn}
\ddot{\bm{y}}_{n+1} = (\tau_n\bm{y}_{n+1}')' + \bm{f}_{n+1} \quad\mbox{in}\quad \mathbb{T}\times(0,T)
\end{equation}
with $\bm{f}_{n+1}=([\partial_s^j,\tau_n]\bm{x}_{n+1}')'$. 
Then, we see that 
\begin{align}\label{EE1}
\frac{\mathrm{d}}{\mathrm{d}t}\|(\dot{\bm{y}}_{n+1},\sqrt{\tau_n}\bm{y}_{n+1}')\|_{L^2}^2
&= (\dot{\tau}_n\bm{y}_{n+1}',\bm{y}_{n+1}')_{L^2} + 2(\dot{\bm{y}}_{n+1},\bm{f}_{n+1})_{L^2} \\
&\leq C_2\|\sqrt{\tau_n}\bm{y}_{n+1}'\|_{L^2}^2 + 2\|\dot{\bm{y}}_{n+1}\|_{L^2}\|\bm{f}_{n+1}\|_{L^2}. \nonumber
\end{align}
Here, we see also that 
\begin{equation}\label{EstLOT}
\|\bm{f}_{n+1}\|_{L^2} \lesssim \|\tau_n'\|_{L^\infty}\|\partial_s^{j+1}\bm{x}_{n+1}\|_{L^2}+\|\partial_s^{j+1}\tau_n\|_{L^2}\|\bm{x}_{n+1}'\|_{L^\infty}.
\end{equation}
Particularly, in the case $j=1,2$, by Poincar\'e's inequality $\|u'\|_{L^2}\leq\|u''\|_{L^2}$ and the Sobolev embedding theorem 
we have $\|\bm{f}_{n+1}\|_{L^2} \lesssim \|\tau_n\|_{H^3}\|\sqrt{\tau_n}\bm{y}_{n+1}'\|_{L^2}$. 
Plugging this into \eqref{EE1} and using Gronwall's inequality, we obtain 
\begin{equation}\label{EE2}
\|(\dot{\bm{x}}_{n+1},\bm{x}_{n+1}')(t)\|_{H^2} \leq C_0\mathrm{e}^{C_2t}\|(\bm{x}_1^\mathrm{in},\bm{x}_0^{\mathrm{in}\,\prime})\|_{H^2}
 \quad\mbox{for}\quad t\in[0,T].
\end{equation}
Moreover, it follows from \eqref{IVPn} that $\|\ddot{\bm{x}}_{n+1}\|_{H^1}\leq\|\tau_n\bm{x}_{n+1}'\|_{H^2} \leq C\|\tau_n\|_{H^2}\|\bm{x}_{n+1}'\|_{H^2}$, 
so that 
\begin{equation}\label{EE3}
\|\ddot{\bm{x}}_{n+1}(t)\|_{H^1} \leq C_1\|\bm{x}_{n+1}'(t)\|_{H^2} \quad\mbox{for}\quad t\in[0,T].
\end{equation}

In view of these estimates \eqref{EE2} and \eqref{EE3}, we define the constants $M_1$ and $M_2$ by $M_1=2C_0M_*$ and $M_2=C_1M_1$, 
and then choosing $T\in(0,T_0]$ so small that $C_2T\leq\log2$ we see that $\bm{x}_{n+1}$ satisfies the estimates in the second line in \eqref{SolEst}. 
Moreover, we have 
\[
\begin{cases}
 \|\bm{x}_{n+1}''(t)-\bm{x}_0^{\mathrm{in}\,\prime\prime}\|_{L^2} \leq \int_0^t\|\dot{\bm{x}}_{n+1}''(t')\|_{L^2}\mathrm{d}t' \leq M_1t, \\
 \|\dot{\bm{x}}_{n+1}'(t)-\bm{x}_1^{\mathrm{in}\,\prime}\|_{L^2} \leq \int_0^t\|\ddot{\bm{x}}_{n+1}'(t')\|_{L^2}\mathrm{d}t' \leq M_2t.
\end{cases}
\]
Here, without loss of generality we may assume that $M_1\leq M_2$. 
Therefore, by further restricting $T$ so small that $2M_0M_2T\leq1$, we see that $\bm{x}_{n+1}$ also satisfies the estimates in the first line in \eqref{SolEst}.

We will show that these approximate solutions also satisfy 
\begin{equation}\label{SolEst2}
\|(\dot{\bm{x}}_n,\bm{x}_n')(t)\|_{H^{m-1}} \leq M_3 \quad\mbox{for}\quad t\in[0,T]
\end{equation}
with a positive constant $M_3$ independent of $n$. 
We recall the $\bm{y}_{n+1}=\partial_s^j\bm{x}_{n+1}$ satisfies \eqref{Eqn} and \eqref{EE1} with \eqref{EstLOT}. 
Now, for $j=1,2,\ldots,m-1$ we have $\|\bm{f}_{n+1}\|_{L^2} \leq C(m,M_0,M_1)(\|\bm{y}_{n+1}'\|_{L^2}+\|\tau_n\|_{H^m})$, 
so that 
\begin{align*}
&\|(\dot{\bm{y}}_{n+1},\sqrt{\tau_n}\bm{y}_{n+1}')(t)\|_{L^2}^2 \\
&\leq \|(\dot{\bm{y}}_{n+1},\sqrt{\tau_n}\bm{y}_{n+1}')(0)\|_{L^2}^2
 + C_3\int_0^t\left( \|(\dot{\bm{y}}_{n+1},\bm{y}_{n+1}')(t')\|_{L^2}^2 + \|\tau_n(t')\|_{H^m}^2 \right)\mathrm{d}t'.
\end{align*}
On the other hand, by Lemma \ref{lem:tau2} we have 
\begin{equation}\label{SolEst3}
\|\tau_n\|_{H^m} \leq C(m,M_0,M_1)\|(\dot{\bm{x}}_n,\bm{x}_n')\|_{H^{m-1}}. 
\end{equation}
Therefore, we obtain 
\begin{align*}
& \|(\dot{\bm{x}}_{n+1},\bm{x}_{n+1}')(t)\|_{H^{m-1}}^2 \\
&\leq C_0\|(\bm{x}_1^\mathrm{in},\bm{x}_0^{\mathrm{in}\,\prime})\|_{H^{m-1}}^2
 + C_3\int_0^t\left( \|(\dot{\bm{x}}_{n+1},\bm{x}_{n+1}')(t')\|_{H^{m-1}}^2 + \|(\dot{\bm{x}}_n,\bm{x}_n')(t')\|_{H^{m-1}}^2 \right)\mathrm{d}t',
\end{align*}
which together with Gronwall's inequality implies 
\begin{align*}
& \|(\dot{\bm{x}}_{n+1},\bm{x}_{n+1}')(t)\|_{H^{m-1}}^2\mathrm{e}^{-C_3t} 
\leq C_0\|(\bm{x}_1^\mathrm{in},\bm{x}_0^{\mathrm{in}\,\prime})\|_{H^{m-1}}^2
 + C_3\int_0^t\|(\dot{\bm{x}}_n,\bm{x}_n')(t')\|_{H^{m-1}}^2 \mathrm{e}^{-C_3t'}\mathrm{d}t'
\end{align*}
for $n=0,1,2,\ldots$. 
Thus, we get 
\begin{align*}
&\|(\dot{\bm{x}}_n,\bm{x}_n')(t)\|_{H^{m-1}}^2\mathrm{e}^{-C_3t} \\
&\leq C_0\|(\bm{x}_1^\mathrm{in},\bm{x}_0^{\mathrm{in}\,\prime})\|_{H^{m-1}}^2\sum_{k=0}^{n-1}\frac{(C_3t)^k}{k!}
 + \frac{(C_3t)^n}{n!}\max_{t'\in[0,t]}\|(\dot{\bm{x}}_0,\bm{x}_0')(t')\|_{H^{m-1}}^2 \\
&\leq C_0\mathrm{e}^{C_3t}\|(\bm{x}_1^\mathrm{in},\bm{x}_0^{\mathrm{in}\,\prime})\|_{H^{m-1}}^2,
\end{align*}
so that $\|(\dot{\bm{x}}_n,\bm{x}_n')(t)\|_{H^{m-1}} \leq C_0\mathrm{e}^{C_3t}\|(\bm{x}_1^\mathrm{in},\bm{x}_0^{\mathrm{in}\,\prime})\|_{H^{m-1}}$ for $t\in[0,T]$. 
Therefore, taking the constant $M_3$ so large that $C_0\mathrm{e}^{C_3T}\|(\bm{x}_1^\mathrm{in},\bm{x}_0^{\mathrm{in}\,\prime})\|_{H^{m-1}}\leq M_3$, 
we obtain the uniform bound \eqref{SolEst2}.

\medskip
{\bf Convergence of the approximate solutions}. 
We proceed to show that the sequence of the approximate solutions $\{(\bm{x}_n,\tau_n)\}$ converges to a solution of \eqref{TFP}--\eqref{IC_bis}. 
Put $\bm{y}_n=\bm{x}_{n+1}-\bm{x}_n$ and $\nu_n=\tau_{n+1}-\tau_n$. 
Then, we see that $\bm{y}_{n+1}$ is a solution of the initial value problem 
\[
\begin{cases}
 \ddot{\bm{y}}_{n+1}=(\tau_{n+1}\bm{y}_{n+1}')' + \bm{f}_{n+1} &\mbox{in}\quad \mathbb{T}\times(0,T), \\
 (\bm{y},\dot{\bm{y}})|_{t=0}=(\bm{0},\bm{0}) &\mbox{in}\quad \mathbb{T}
\end{cases}
\]
with $\bm{f}_{n+1}=(\nu_n\bm{x}_{n+1}')'$. 
Here, we have $\|\bm{f}_{n+1}\|_{L^2} \leq C\|\nu_n'\|_{H^1}\|\bm{x}_{n+1}'\|_{H^1} \leq C_0\|\nu_n'\|_{H^1}$. 
Therefore, carrying out the same energy estimate as in \eqref{EE1}, we obtain 
\[
\|(\dot{\bm{y}}_{n+1},\bm{y}_{n+1}')(t)\|_{L^2}^2 \leq C_2\int_0^t\|\nu_n(t')\|_{H^1}^2\mathrm{d}t'.
\]
On the other hand, $\nu_n$ is a solution of the periodic boundary value problem 
\begin{align*}
-\nu_n'' + |\bm{x}_{n+1}|^2\nu_n = h_n \quad\mbox{in}\quad \mathbb{T}\times(0,T)
\end{align*}
with 
\begin{align*}
h_n
&= (\dot{\bm{x}}_{n+1}'+\dot{\bm{x}}_n')\cdot\dot{\bm{y}}_n' -\tau_n(\bm{x}_{n+1}''+\bm{x}_n'')\cdot\bm{y}_n'' \\
&= \left( (\dot{\bm{x}}_{n+1}'+\dot{\bm{x}}_n')\cdot\dot{\bm{y}}_n -\tau_n(\bm{x}_{n+1}''+\bm{x}_n'')\cdot\bm{y}_n' \right)' \\
&\quad\;
  -(\dot{\bm{x}}_{n+1}'+\dot{\bm{x}}_n')'\cdot\dot{\bm{y}}_n + \left(\tau_n(\bm{x}_{n+1}''+\bm{x}_n'')\right)'\cdot\bm{y}_n'.
\end{align*}
By Lemmas \ref{lem:BVP1} and \ref{lem:BVP2}, we obtain 
\begin{align}\label{Est_nu}
\|\nu_n\|_{H^1}
&\leq C_0\bigl( \|(\dot{\bm{x}}_{n+1}'+\dot{\bm{x}}_n')\cdot\dot{\bm{y}}_n -\tau_n(\bm{x}_{n+1}''+\bm{x}_n'')\cdot\bm{y}_n'\|_{L^2} \\
&\quad\;
 + \|(\dot{\bm{x}}_{n+1}'+\dot{\bm{x}}_n')'\cdot\dot{\bm{y}}_n - \left(\tau_n(\bm{x}_{n+1}''+\bm{x}_n'')\right)'\cdot\bm{y}_n'\|_{L^1} \bigr) \nonumber \\
&\leq C_0( \|\dot{\bm{x}}_{n+1}'+\dot{\bm{x}}_n'\|_{H^1}\|\dot{\bm{y}}_n\|_{L^2}
 + \|\tau_n(\bm{x}_{n+1}''+\bm{x}_n'')\|_{H^1}\|\bm{y}_n'\|_{L^2} ) \nonumber \\
&\leq C_1\|(\dot{\bm{y}}_n,\bm{y}_n')\|_{L^2}, \nonumber
\end{align}
so that $\|(\dot{\bm{y}}_{n+1},\bm{y}_{n+1}')(t)\|_{L^2}^2 \leq C_2\int_0^t\|(\dot{\bm{y}}_n,\bm{y}_n')(t)\|_{L^2}^2\mathrm{d}t'$. 
Therefore, we get 
\begin{equation}\label{SolEst4}
\max_{t\in[0,T]}\|(\dot{\bm{y}}_n,\bm{y}_n')(t)\|_{L^2}^2 \leq \frac{(C_2T)^n}{n!}\max_{t\in[0,T]}\|(\dot{\bm{y}}_0,\bm{y}_0')(t)\|_{L^2}^2 
 \quad\mbox{for}\quad n=0,1,2,\ldots.
\end{equation}
This shows that $\{(\dot{\bm{x}}_n,\bm{x}_n')\}$ is a Cauchy sequence in $C^0([0,T];L^2)$. 
This together with the uniform bound \eqref{SolEst2} and an interpolation inequality implies that $\{\bm{x}_n\}$ is a Cauchy sequence in 
$C^0([0,T];H^{m-1})\cap C^1([0,T];H^{m-2})$. 
We denote by $\bm{x}$ its limit. 
By standard compactness arguments, we have also $\bm{x}\in L^\infty(0,T;H^m)$ and $\dot{\bm{x}}\in L^\infty(0,T;H^{m-1})$. 
In view of \eqref{Est_nu} and \eqref{SolEst4}, we see also that $\{\tau_n\}$ is a Cauchy sequence in $C^0([0,T];H^1)$. 
This together with the uniform bound \eqref{SolEst3} implies that $\{\tau_n\}$ is a Cauchy sequence in $C^0([0,T];H^{m-1})$. 
We denote its limit by $\tau$, which is also in $L^\infty(0,T;H^m)$ and satisfies $\tau(s,t) \geq C_0^{-1}$ for $(s,t)\in\mathbb{T}\times[0,T]$. 
Moreover, by the uniform bound $\|\dot{\tau}_n\|_{L^\infty(\mathbb{T}\times(0,T))} \leq C_2$, we have also $\dot{\tau}\in L^\infty(\mathbb{T}\times(0,T))$. 
Obviously, $(\bm{x},\tau)$ is a solution of the initial value problem \eqref{TFP}--\eqref{IC_bis}. 
Since we can regard $\bm{x}$ as a unique solution of the initial value problem to linear hyperbolic equation \eqref{TFP} and \eqref{IC_bis}, 
we have enough regularity on the coefficient $\tau$ to ensure that $\bm{x}\in C^0([0,T];H^m)\cap C^1([0,T];H^{m-1})$. 
Then, since $\tau$ is a solution of the periodic boundary value problem \eqref{TFBVP}, we get $\tau\in C^0([0,T];H^m)$. 
Finally, the uniqueness of the solution can be proved in the same way as above. 
\end{proof}

Now, Theorem \ref{Th:LWP} follows directly from Proposition \ref{Prop:WP}. 
In fact, since the initial data are supposed to satisfy \eqref{Compatibility_Condition}, 
by Lemma \ref{lem:Fenchel} we have $\|\bm{x}_0^{\mathrm{in}\,\prime\prime}\|_{L^2} \geq \sqrt{2\pi}$, 
so that we can apply Proposition \ref{Prop:WP} to ensure the existence of the solution $(\bm{x},\tau)$ to \eqref{TFP}--\eqref{IC_bis}. 
Then, by Proposition \ref{Prop:Equivalence} we see that $(\bm{x},\tau)$ is a unique solution of the original problem \eqref{Eq}--\eqref{IC}.

\section{Blow-up criterion}\label{sec:BUC}
In this last section, we will prove Theorem \ref{Th:BUC}. 
As was explained in Section \ref{sec:intro}, if the initial velocity $\bm{x}_1^\mathrm{in}(s)$ is constant, 
then the initial value problem \eqref{Eq}--\eqref{IC} has a unique time global solution. 
Therefore, it is sufficient to consider the case where $\bm{x}_1^\mathrm{in}(s)$ is not constant, 
that is, the case $\|P_1\bm{x}_1^\mathrm{in}\|_{L^2}>0$, 
where $P_1$ is the orthogonal projection defined in Lemma \ref{lem:Below}. 
Since we have already obtained the local well-posedness of the initial value problem \eqref{Eq}--\eqref{IC} in Theorem \ref{Th:LWP}, 
it is sufficient to show the following proposition.

\begin{proposition}\label{Prop:BUC}
For any integer $m\geq3$ and positive constants $M$ and $T$, there exists a positive constant $C$ such that the following holds true: 
Suppose that the initial data $(\bm{x}_0^\mathrm{in},\bm{x}_1^\mathrm{in})$ satisfy 
\[
\begin{cases}
 \|\bm{x}_0^\mathrm{in}\|_{H^m} + \|\bm{x}_1^\mathrm{in}\|_{H^{m-1}} \leq M, \\
 |\bm{x}_0^{\mathrm{in}\,\prime}(s)|\equiv1, \quad \|P_1\bm{x}_1^\mathrm{in}\|_{L^2}\geq M^{-1}.
\end{cases}
\]
Then, any regular solution $(\bm{x},\tau)$ of the initial value problem \eqref{Eq}--\eqref{IC} satisfying 
\[
\int_0^T \|\bm{x}''(t)\|_{H^1}^2\mathrm{d}t \leq M
\]
satisfies $\|\bm{x}(t)\|_{H^m}+\|\dot{\bm{x}}(t)\|_{H^{m-1}} \leq C$ for $t\in[0,T]$. 
\end{proposition}

To show this proposition, the following lemma is crucial.

\begin{lemma}\label{lem:BUC}
There exists an absolute positive constant $C$ such that any regular solution $(\bm{x},\tau)$ of the initial value problem \eqref{Eq}--\eqref{IC} satisfies 
\begin{align*}
&\|\dot{\bm{x}}''(t)\|_{L^2}^2 + \|\sqrt{\tau(t)}\bm{x}'''(t)\|_{L^2}^2 + 4\||\dot{\bm{x}}'(t)||\bm{x}''(t)|\|_{L^2}^2 \\
&\leq \bigl( \|\bm{x}_1^{\mathrm{in}\,\prime\prime}\|_{L^2}^2 + \|\sqrt{\tau_0^\mathrm{in}}\bm{x}_0^{\mathrm{in}\,\prime\prime\prime}\|_{L^2}^2
 + 4\||\bm{x}_1^{\mathrm{in}\,\prime}||\bm{x}_0^{\mathrm{in}\,\prime\prime}|\|_{L^2}^2 \bigr)
 \exp\left( C\|\bm{x}_1^\mathrm{in}\|_{L^2}\int_0^t\|\bm{x}'''(t')\|_{L^2}^2\mathrm{d}t' \right),
\end{align*}
where $\tau_0^\mathrm{in}$ is the initial tension of the string determined from $(\bm{x}_0^\mathrm{in},\bm{x}_1^\mathrm{in})$. 
\end{lemma}

\begin{proof}
We first recall that the tension $\tau$ satisfies necessarily the periodic boundary value problem \eqref{BVP}. 
By the equation of motion \eqref{Eq} and integration by parts, we see that 
\begin{align*}
\frac{\mathrm{d}}{\mathrm{d}t}( \|\dot{\bm{x}}''\|_{L^2}^2 + \|\sqrt{\tau}\bm{x}'''\|_{L^2}^2 )
&= 2(\dot{\bm{x}}'',\ddot{\bm{x}}'')_{L^2} + 2 (\tau\bm{x}''',\dot{\bm{x}}''')_{L^2} + (\dot{\tau}\bm{x}''',\bm{x}''')_{L^2} \\
&= 2(\dot{\bm{x}}'',(\tau\bm{x}')'''-(\tau\bm{x}''')')_{L^2} + (\dot{\tau}\bm{x}''',\bm{x}''')_{L^2} \\
&= 2(\dot{\bm{x}}'',2\tau'\bm{x}'''+3\tau''\bm{x}''+\tau'''\bm{x}')_{L^2} + (\dot{\tau}\bm{x}''',\bm{x}''')_{L^2}.
\end{align*}
By differentiating the constraint \eqref{Constraint} with respect to $s$ and $t$, we obtain $\bm{x}'\cdot\dot{\bm{x}}''+\dot{\bm{x}}'\cdot\bm{x}''=0$, 
so that $(\dot{\bm{x}}'',\tau'''\bm{x}')_{L^2} = -(\dot{\bm{x}}'\cdot\bm{x}'',\tau''')_{L^2}
 = (\dot{\bm{x}}''\cdot\bm{x}''+\dot{\bm{x}}'\cdot\bm{x}''',\tau'')_{L^2}$. 
Plugging this into the above equation and using \eqref{BVP} to eliminate $\tau''$, we get 
\begin{align*}
&\frac{\mathrm{d}}{\mathrm{d}t}( \|\dot{\bm{x}}''\|_{L^2}^2 + \|\sqrt{\tau}\bm{x}'''\|_{L^2}^2 ) \\
&= 4(\dot{\bm{x}}'',\tau'\bm{x}''')_{L^2} + 2(4\dot{\bm{x}}''\cdot\bm{x}''+\dot{\bm{x}}'\cdot\bm{x}''',\tau|\bm{x}''|^2-|\dot{\bm{x}}'|^2)_{L^2}
 + (\dot{\tau}\bm{x}''',\bm{x}''')_{L^2}.
\end{align*}
On the other hand, by \eqref{Eq} we see also that 
\begin{align*}
4\frac{\mathrm{d}}{\mathrm{d}t}\||\dot{\bm{x}}'||\bm{x}''|\|_{L^2}^2
&= 8(\dot{\bm{x}}'\cdot\ddot{\bm{x}}',|\bm{x}''|^2)_{L^2} + 8(\bm{x}''\cdot\dot{\bm{x}}'',|\dot{\bm{x}}'|^2)_{L^2} \\
&= 8(\dot{\bm{x}}'\cdot(\tau\bm{x}')'',|\bm{x}''|^2)_{L^2} + 8(\bm{x}''\cdot\dot{\bm{x}}'',|\dot{\bm{x}}'|^2)_{L^2} \\
&= 8(\dot{\bm{x}}'\cdot(\tau\bm{x}'''+2\tau'\bm{x}''),|\bm{x}''|^2)_{L^2} + 8(\bm{x}''\cdot\dot{\bm{x}}'',|\dot{\bm{x}}'|^2)_{L^2},
\end{align*}
where we used the identity $\dot{\bm{x}}'\cdot\bm{x}'=0$. 
By adding the above two identities, the troubling term $(\bm{x}''\cdot\dot{\bm{x}}'',|\dot{\bm{x}}'|^2)_{L^2}$ is cancelled and we obtain 
\begin{align}\label{EI}
& \frac{\mathrm{d}}{\mathrm{d}t}( \|\dot{\bm{x}}''\|_{L^2}^2 + \|\sqrt{\tau}\bm{x}'''\|_{L^2}^2 + 4\||\dot{\bm{x}}'||\bm{x}''|\|_{L^2}^2) \\
&= 4(\dot{\bm{x}}'',\tau'\bm{x}''')_{L^2} + 8(\dot{\bm{x}}'',\tau|\bm{x}''|^2\bm{x}'')_{L^2} + 10(\dot{\bm{x}}',\tau|\bm{x}''|^2\bm{x}''')_{L^2} \nonumber \\
&\quad\;
 + 16(\dot{\bm{x}}',\tau'|\bm{x}''|^2\bm{x}'')_{L^2} - 2(|\dot{\bm{x}}'|^2\dot{\bm{x}}',\bm{x}''')_{L^2} + (\dot{\tau}\bm{x}''',\bm{x}''')_{L^2}. \nonumber
\end{align}

Before evaluating each term in the right-hand side of the above identity, we prepare several estimates. 
Recall that we have the conservation of the kinetic energy $\|\dot{\bm{x}}(t)\|_{L^2}\equiv\|\bm{x}_1^{\mathrm{in}}\|_{L^2}$. 
By integration by parts, we see also that $\|\dot{\bm{x}}'(t)\|_{L^2}^2=-(\dot{\bm{x}}(t),\dot{\bm{x}}''(t))_{L^2}\leq
 \|\dot{\bm{x}}(t)\|_{L^2}\|\dot{\bm{x}}''(t)\|_{L^2}=\|\bm{x}_1^{\mathrm{in}}\|_{L^2}\|\dot{\bm{x}}''(t)\|_{L^2}$. 
By Lemma \ref{lem:Fenchel} we have $\|\bm{x}''(t)\|_{L^2}\geq\sqrt{2\pi}$. 
Therefore, by Lemma \ref{lem:tau1} we obtain $\|\tau(t)\|_{W^{1,\infty}}\lesssim\|\dot{\bm{x}}'(t)\|_{L^2}^2$. 
To evaluate $\dot{\tau}$, we see that it satisfies $-\dot{\tau}''+|\bm{x}''|^2\dot{\tau}=h$ with 
\begin{align*}
h &= 2\dot{\bm{x}}'\cdot\ddot{\bm{x}}'-2\tau\bm{x}''\cdot\dot{\bm{x}}'' \\
&= 2\dot{\bm{x}}'\cdot(\tau\bm{x}')''-2\tau\bm{x}''\cdot\dot{\bm{x}}'' \\
&= 2\dot{\bm{x}}'\cdot(\tau\bm{x}'''+2\tau'\bm{x}'')-2\tau\bm{x}''\cdot\dot{\bm{x}}'' \\
&= (4\tau\dot{\bm{x}}'\cdot\bm{x}'')' - 2(\tau\dot{\bm{x}}'\cdot\bm{x}'''+3\tau\bm{x}''\cdot\dot{\bm{x}}''),
\end{align*}
where we used the identity $\dot{\bm{x}}'\cdot\bm{x}'=0$. 
Therefore, by Lemmas \ref{lem:BVP1}--\ref{lem:BVP2} we obtain 
\begin{align*}
\|\dot{\tau}\|_{L^\infty}
&\lesssim \|\tau\dot{\bm{x}}'\cdot\bm{x}''\|_{L^1} + \|\tau\dot{\bm{x}}'\cdot\bm{x}'''\|_{L^1} + \|\tau\bm{x}''\cdot\dot{\bm{x}}''\|_{L^1} \\
&\lesssim \|\sqrt{\tau}\|_{L^\infty}( \|\sqrt{\tau}\bm{x}''\|_{L^2}\|\dot{\bm{x}}'\|_{H^1} + \|\sqrt{\tau}\bm{x}'''\|_{L^2}\|\dot{\bm{x}}'\|_{L^2} ).
\end{align*}
Here, by integrating \eqref{BVP} with respect to $s$ over $\mathbb{T}$, we have $\|\sqrt{\tau}\bm{x}''\|_{L^2}=\|\dot{\bm{x}}'\|_{L^2}$. 
Note also Poincar\'e's inequality $\|u'\|_{L^2}\leq\|u''\|_{L^2}$. 
Therefore, we obtain 
\begin{align*}
\|\dot{\tau}\|_{L^\infty}
&\lesssim \|\dot{\bm{x}}'\|_{L^2}^2 ( \|\dot{\bm{x}}''\|_{L^2} + \|\sqrt{\tau}\bm{x}'''\|_{L^2} ) \\
&\lesssim \|\bm{x}_1^\mathrm{in}\|_{L^2} \|\dot{\bm{x}}''\|_{L^2} ( \|\dot{\bm{x}}''\|_{L^2} + \|\sqrt{\tau}\bm{x}'''\|_{L^2} ).
\end{align*}

Now, we go back to evaluate each term in the right-hand side of \eqref{EI}. 
In the following, we will use the above estimates and the Gagliardo--Nirenberg inequality 
$\|u'\|_{L^6}\lesssim\|u''\|_{L^2}^{\frac23}\|u\|_{L^2}^{\frac13}$ without any comment. 
\begin{align*}
|(\dot{\bm{x}}'',\tau'\bm{x}''')_{L^2}|
&\leq \|\dot{\bm{x}}''\|_{L^2}\|\tau'\|_{L^\infty}\|\bm{x}'''\|_{L^2} \\
&\leq \|\bm{x}_1^{\mathrm{in}}\|_{L^2}\|\bm{x}'''\|_{L^2}\|\dot{\bm{x}}''\|_{L^2}^2,
\end{align*}
\begin{align*}
|(\dot{\bm{x}}'',\tau|\bm{x}''|^2\bm{x}'')_{L^2}|
&\leq \|\dot{\bm{x}}''\|_{L^2}\|\sqrt{\tau}\|_{L^\infty}\|\sqrt{\tau}\bm{x}''\|_{L^2}\|\bm{x}''\|_{L^\infty}^2 \\
&\lesssim \|\dot{\bm{x}}'\|_{L^2}^2\|\bm{x}''\|_{H^1}^2\|\dot{\bm{x}}''\|_{L^2} \\
&\leq \|\bm{x}_1^\mathrm{in}\|_{L^2}\|\bm{x}''\|_{H^1}^2\|\dot{\bm{x}}''\|_{L^2}^2,
\end{align*}
\begin{align*}
|(\dot{\bm{x}}',\tau|\bm{x}''|^2\bm{x}''')_{L^2}|
&\leq \|\dot{\bm{x}}'\|_{L^2}\|\sqrt{\tau}\|_{L^\infty}\|\bm{x}''\|_{L^\infty}^2\|\sqrt{\tau}\bm{x}'''\|_{L^2} \\
&\lesssim \|\dot{\bm{x}}'\|_{L^2}^2\|\bm{x}''\|_{H^1}^2\|\sqrt{\tau}\bm{x}'''\|_{L^2} \\
&\leq \|\bm{x}_1^\mathrm{in}\|_{L^2}\|\bm{x}''\|_{H^1}^2\|\dot{\bm{x}}''\|_{L^2}\|\sqrt{\tau}\bm{x}'''\|_{L^2},
\end{align*}
\begin{align*}
|(\dot{\bm{x}}',\tau'|\bm{x}''|^2\bm{x}'')_{L^2}|
&\leq \|\dot{\bm{x}}'\|_{L^2}\|\tau'\|_{L^\infty}\|\bm{x}''\|_{L^6}^3 \\
&\lesssim \|\dot{\bm{x}}'\|_{L^2}^3\|\bm{x}'\|_{L^2}\|\bm{x}'''\|_{L^2}^2 \\
&\lesssim \|\bm{x}_1^\mathrm{in}\|_{L^2}\|\bm{x}'''\|_{L^2}^2\|\dot{\bm{x}}''\|_{L^2}^2,
\end{align*}
\begin{align*}
|(|\dot{\bm{x}}'|^2\dot{\bm{x}}',\bm{x}''')_{L^2}|
&\leq \|\dot{\bm{x}}'\|_{L^6}^3\|\bm{x}'''\|_{L^2} \\
&\lesssim \|\dot{\bm{x}}\|_{L^2}\|\dot{\bm{x}}''\|_{L^2}^2\|\bm{x}'''\|_{L^2} \\
&= \|\bm{x}_1^\mathrm{in}\|_{L^2}\|\bm{x}'''\|_{L^2}\|\dot{\bm{x}}''\|_{L^2}^2,
\end{align*}
\begin{align*}
|(\dot{\tau}\bm{x}''',\bm{x}''')_{L^2}|
&\leq \|\dot{\tau}\|_{L^\infty}\|\bm{x}'''\|_{L^2}^2 \\
&\lesssim \|\bm{x}_1^\mathrm{in}\|_{L^2}\|\bm{x}'''\|_{L^2}^2( \|\dot{\bm{x}}''\|_{L^2}^2 + \|\sqrt{\tau}\bm{x}'''\|_{L^2}^2 ).
\end{align*}
We note also that $\sqrt{2\pi}\leq\|\bm{x}''\|_{L^2}\leq\|\bm{x}'''\|_{L^2}$. 
Plugging the above estimates into \eqref{EI}, we obtain 
\begin{align*}
& \frac{\mathrm{d}}{\mathrm{d}t}( \|\dot{\bm{x}}''\|_{L^2}^2 + \|\sqrt{\tau}\bm{x}'''\|_{L^2}^2 + 4\||\dot{\bm{x}}'||\bm{x}''|\|_{L^2}^2) 
\lesssim \|\bm{x}_1^\mathrm{in}\|_{L^2}\|\bm{x}'''\|_{L^2}^2( \|\dot{\bm{x}}''\|_{L^2}^2 + \|\sqrt{\tau}\bm{x}'''\|_{L^2}^2 ),
\end{align*}
which yields the desired estimate. 
\end{proof}

\begin{proof}[Proof of Proposition \ref{Prop:BUC}]
It is sufficient to show this proposition in the case $m=3$. 
In fact, once we obtain the bound $\|\bm{x}(t)\|_{H^3}+\|\dot{\bm{x}}(t)\|_{H^2}\leq C$ for $t\in[0,T]$, 
the higher order bound $\|\bm{x}(t)\|_{H^m}+\|\dot{\bm{x}}(t)\|_{H^{m-1}}\leq C$ for $t\in[0,T]$ can follow in almost the same way as the proof 
of the uniform bound \eqref{SolEst2} of the approximate solutions carried out in the proof of Proposition \ref{Prop:WP}. 
Therefore, we consider the case $m=3$. 
In the following, we simply denote $C=C(M,T)$, which may change from line to line. 
Recall that we have $\|\bm{x}''(t)\|_{L^2}\geq\sqrt{2\pi}$ and $\|\dot{\bm{x}}(t)\|_{L^2}=\|\bm{x}_1^\mathrm{in}\|_{L^2}$. 
By Lemma \ref{lem:tau1}, we have $\|\tau_0^\mathrm{in}\|_{L^\infty} \lesssim \|\bm{x}_1^{\mathrm{in}\,\prime}\|_{L^2}^2$. 
By Sobolev embedding theorem, we have 
$\||\bm{x}_1^{\mathrm{in}\,\prime}||\bm{x}_0^{\mathrm{in}\,\prime\prime}|\|_{L^2}
 \lesssim \|\bm{x}_1^{\mathrm{in}\,\prime}\|_{H^1}\|\bm{x}_0^{\mathrm{in}\,\prime\prime}\|_{H^1}$. 
Therefore, by Lemma \ref{lem:BUC} we obtain $\|\dot{\bm{x}}''(t)\|_{L^2} + \|\sqrt{\tau(t)}\bm{x}'''(t)\|_{L^2} \leq C$, 
which implies immediately $\|\dot{\bm{x}}(t)\|_{H^2}\leq C$, so that 
$\|\bm{x}(t)\|_{H^2}\leq\|\bm{x}_0^\mathrm{in}\|_{H^2}+\int_0^t\|\dot{\bm{x}}(t')\|_{H^2}\mathrm{d}t' \leq C$. 
Moreover, by Lemmas \ref{lem:tau1} and \ref{lem:Below} we get $\min_{s\in\mathbb{T}}\tau(s,t)\gtrsim\|P_1\bm{x}_1^\mathrm{in}\|_{L^2}^2\geq M^{-2}$, 
which implies that $\|\bm{x}'''(t)\|_{L^2}\leq C$. 
Therefore, we finally obtain $\|\bm{x}(t)\|_{H^3}+\|\dot{\bm{x}}(t)\|_{H^2} \leq C$ for $t\in[0,T]$. 
This completes the proof. 
\end{proof}

\medskip
\noindent
{\bf Acknowledgement} \\
T. I. is partially supported by JSPS KAKENHI Grant Number JP23K22404.


\bigskip
Tatsuo Iguchi \par
{\sc Department of Mathematics} \par
{\sc Faculty of Science and Technology, Keio University} \par
{\sc 3-14-1 Hiyoshi, Kohoku-ku, Yokohama, 223-8522, Japan} \par
E-mail: \texttt{iguchi@math.keio.ac.jp}

\bigskip
Masahiro Takayama \par
{\sc Department of Mathematics} \par
{\sc Faculty of Science and Technology, Keio University} \par
{\sc 3-14-1 Hiyoshi, Kohoku-ku, Yokohama, 223-8522, Japan} \par
E-mail: \texttt{masahiro@math.keio.ac.jp}

\end{document}